\documentclass[12pt,a4paper]{amsart}

\usepackage[utf8]{inputenc}
\usepackage[margin=1in]{geometry}
\usepackage{thmtools} 
\usepackage{hyperref}
\usepackage{etoolbox} 
\usepackage[capitalize,noabbrev]{cleveref}
\usepackage{amsmath}
\usepackage{amsfonts}
\usepackage{amsthm}
\usepackage{amssymb}
\usepackage{xcolor}
\usepackage{graphicx}
\usepackage{tikz-cd} 
\usepackage{bbm} 
\usepackage{enumerate} 
\usepackage{stmaryrd} 
\usepackage{subcaption} 
\usepackage{microtype} 
\usepackage{mathrsfs} 
\usepackage{enumitem}
\usepackage{comment}
\usetikzlibrary{decorations.pathreplacing}

\theoremstyle{plain}
\newtheorem{thrm}{Theorem}[section]
\def\bthm{\begin{thrm}}
\def\ethm{\end{thrm}}

\newtheorem{theoremx}{Theorem}

\newtheorem{prop}[thrm]{Proposition}
\def\bprop{\begin{prop}}
\def\eprop{\end{prop}}

\newtheorem{ques}[thrm]{Question}
\def\bques{\begin{ques}}
\def\eques{\end{ques}}

\newtheorem{cjec}[thrm]{Conjecture}
\def\bcjec{\begin{cjec}}
\def\ecjec{\end{cjec}}

\newtheorem{cor}[thrm]{Corollary}
\def\bcor{\begin{cor}}
\def\ecor{\end{cor}}

\newtheorem{fact}[thrm]{Fact}
\def\bfact{\begin{fact}}
\def\efact{\end{fact}}

\newtheorem{lem}[thrm]{Lemma}
\def\blem{\begin{lem}}
\def\elem{\end{lem}}

\theoremstyle{definition}
\newtheorem{defn}[thrm]{Definition}
\def\bdefn{\begin{defn}}
\def\edefn{\end{defn}}

\newtheorem*{conc}{Conclusion}
\def\bconc{\begin{conc}}
\def\econc{\end{conc}}

\newtheorem{alg}[thrm]{Algorithm}
\def\balg{\begin{alg}}
\def\ealg{\end{alg}}

\def\bproof{\begin{proof}}
\def\eproof{\end{proof}}

\theoremstyle{remark}

\newtheorem{rem}[thrm]{Remark}
\def\brem{\begin{rem}}
\def\erem{\end{rem}}

\newtheorem{ex}[thrm]{Example}
\def\bex{\begin{ex}}
\def\eex{\end{ex}}

\newtheorem{exs}[thrm]{Examples}
\def\bexs{\begin{exs}}
\def\eexs{\end{exs}}

\newtheorem{obs}{Observation}
\def\bobs{\begin{obs}}
\def\eobs{\end{obs}}	

\let\ra=\rightarrow

\let\a=\alpha

\let\d=\partial
\let\setminus=\smallsetminus

\def\A{\mathcal{A}}

\def\G{\mathcal{G}}

\def\N{\mathbb{N}}
\def\Z{\mathbb{Z}}
\def\R{\mathbb{R}}

\def\K{\mathbb{K}}

\def\l{\mathds{l}}

\def\d{\partial}

\def\p{\mathfrak{p}}

\DeclareMathOperator{\id}{id}

\DeclareMathOperator{\Supp}{Supp}

\DeclareMathOperator{\init}{in}

\DeclareMathOperator{\reg}{reg}
\DeclareMathOperator{\Ind}{Ind}
\DeclareMathOperator{\ind}{Ind}

\DeclareMathOperator{\cdeg}{cdeg}
\DeclareMathOperator{\Cut}{Cut}
\DeclareMathOperator{\lex}{lex}

\newcommand{\card}[1]{\lvert #1 \rvert}

\newcommand{\ul}[1]{\underline{#1}}
\newcommand{\mc}[1]{\mathcal{#1}}

\newcommand{\wh}[1]{\widehat{#1}}

\renewcommand{\subset}{\subseteq}

\AtBeginDocument{%
   \def\MR#1{}
}

\newcommand{\edgelabeledgraph}[2]{
	\begin{tikzpicture}[scale = 1.2]          
	\newcommand*\points{#1}     
	\newcommand*\edges{#2}          
	\newcommand*\scale{0.75}          
	\foreach \x/\y/\z/\w in \points {
		\draw[fill = black!50] (\scale*\x,\scale*\y) circle [radius = 0.1] node[label = {[label distance = 0.05 cm]\w: $\z$}] (\z) {}; 
	}
	\foreach \x/\y/\p/\a/\l in \edges { \ifthenelse{\x = \y}{\draw (\x) to[out = \a -20, in = \a + 20, looseness = 50]  node [pos = \p, label =\a: \textcolor{blue}{\l}] {} (\y); }{
    \draw (\x) -- (\y) node [pos = \p, label =\a: \textcolor{blue}{\l}] {};}
    }
	\end{tikzpicture}
}

\begin{document}

\title[Castelnuovo--Mumford Regularity of Binomial Edge Ideals]{Combinatorics of Castelnuovo--Mumford Regularity of Binomial Edge Ideals}
\author{Adam LaClair}
\thanks{The author was partially supported by NSF grant DMS-2100288 and by Simons Foundation Collaboration Grant for Mathematicians \#580839}

\maketitle

\begin{abstract}
Since the introduction of binomial edge ideals $J_{G}$ by Herzog et al.\ and independently Ohtani, there has been significant interest in relating algebraic invariants of the binomial edge ideal with combinatorial invariants of the underlying graph $G$. Here, we take up a question considered by Herzog and Rinaldo regarding Castelnuovo--Mumford regularity of block graphs. To this end, we introduce a new invariant $\nu(G)$ associated to any simple graph $G$, defined as the  maximal total length of a certain collection of induced paths within $G$ subject to conditions on the induced subgraph. We prove that for any graph $G$, $\nu(G) \leq \reg(J_{G})-1$, and that the length of a longest induced path of $G$ is less than or equal to $\nu(G)$; this refines an inequality of Matsuda and Murai. We then investigate the question: when is $\nu(G) = \reg(J_{G})-1$? We prove that equality holds when $G$ is closed; this gives a new characterization of a result of Ene and Zarojanu, and when $G$ is bipartite and $J_{G}$ is Cohen-Macaulay; this gives a new characterization of a result of Jayanathan and Kumar. For a block graph $G$, we prove that $\nu(G)$ admits a combinatorial characterization independent of any auxiliary choices, and we prove that $\nu(G) = \reg(J_{G})-1$. This gives $\reg(J_{G})$ a combinatorial interpretation for block graphs, and thus answers the question of Herzog and Rinaldo.
\end{abstract}

\section{Introduction}

Castelnuovo--Mumford regularity, introduced by David Mumford in the 1960s \cite{mumford1966lectures}, is a fundamental invariant in commutative algebra and algebraic geometry that roughly measures how \textit{complicated} a module or sheaf is. It is an interesting and difficult question to provide a combinatorial interpretation of the Castelnuovo--Mumford regularity for families of ideals possessing an underlying combinatorial structure. One such family of ideals, studied extensively over the past decade, is the class of binomial edge ideals $J_{G}$ associated to a graph $G$. These were introduced by Herzog et al.\ \cite{herzog2010binomial} and independently Ohtani \cite{ohtani2011graphs}; see Section \ref{sec_background} for precise definitions. There are elegant combinatorial upper bounds for $\reg(J_{G})$, the Castelnuovo--Mumford regularity of the binomial edge ideal, in terms of: maximum number of clique disjoint edges of the graph \cite{malayeri2021proof} and number of vertices of the graph \cite{kiani2016castelnuovo}. On the other hand, Matsuda and Murai proved that the length of a longest induced path of $G$ gives a lower bound for $\reg(J_{G})-1$ {\cite[Corollary 2.3]{matsuda2013regularity}}. More recently, \cite{ambhore2023mathrm} and \cite{jaramillo2023connected} have investigated the question of giving a lower bound on $\reg(J_{G})-1$ via the $\mathrm{v}$-domination number of binomial edge ideals. Inspired by Matsuda and Murai's lower bound, we ask the question:
\bques
\label{ques_intro}
When is it the case that $\sum_{i=1}^{\ell} \card{E(P_{i})} \leq \reg(J_{G})-1$ for vertex-disjoint induced paths $P_{1},\ldots,P_{\ell}$ of $G$?
\eques
Notice that the case of $\ell=1$ is the Matsuda--Murai lower bound. In general, an arbitrary choice of vertex-disjoint induced paths will not realize a lower bound for $\reg(J_{G})-1$. For instance, the complete graph on $n \geq 2$ vertices has $\reg(J_{G}) = 2$ yet $\lfloor \frac{n}{2} \rfloor$ vertex-disjoint induced edges. Our observations are that (1) it is possible to label the vertices of $G$ so that each path $P_{i}$ corresponds to a monomial in the generating set of $\init_{\lex}(J_{G})$ via Herzog et al's characterization of $\init_{\lex}(J_{G})$; and (2) that with restrictions on the edges appearing between the $P_{i}$'s the corresponding monomials realize a regular sequence whose free resolution is a subcomplex of the free resolution of $\init_{\lex}(J_{G})$. Thus, this provides a lower bound on $\reg(\init_{\lex}(J_{G}))$ in terms of the total number of edges appearing in the $P_{i}$. From this lower bound on the Castelnuovo--Mumford regularity of the initial ideal, we obtain a lower bound on the Castelnuovo--Mumford regularity for $J_{G}$ via Conca--Varbaro's theorem on the preservation of extremal Betti numbers for ideals with squarefree initial ideal. In Section \ref{sec_the_invariant_nu}, we define the invariant $\nu(G)$ along these lines, and we prove the inequality
\begin{align*}
    \nu(G) \leq \reg(J_{G})-1,
\end{align*}
(Theorem \ref{thm_nu_leq_reg}). 

In the remainder of this paper, we take up the question of equality of $\nu(G)$ and $\reg(J_{G})-1$. Various authors have considered the question of describing $\reg(J_{G})$, a purely algebraic invariant, in terms of combinatorial properties of $G$. Ene and Zarojanu showed that $\reg(J_{G})-1$ agrees with the length of the longest induced path of the graph when $G$ is a closed graph \cite{ene2015regularity}. Jayanthan and Kumar gave a combinatorial interpretation of $\reg(J_{G})$ when $G$ is bipartite and $J_{G}$ is Cohen--Macaulay \cite{jayanthan2019regularityCM&BipartiteGraphs}. The graphs having $\reg(J_{G}) \leq 3$ have been classified by Kiani and Saeedi Madani in \cite{kiani2012binomial} and \cite{madani2018binomial}. For the computation of $\reg(J_{G})$ in further cases, see the survey article \cite{mukherjee2022homological}. In Section \ref{sec_equality_nu_and_reg}, we show that $\nu(G)$ agrees with $\reg(J_{G}) -1$ when $G$ is a closed graph (Corollary \ref{cor_nu_equals_reg_closed_graphs}), and when $G$ is bipartite and $J_{G}$ is Cohen-Macaulay (Corollary \ref{cor_nu_equals_reg_CM_bipartite_graphs}). 

In Sections \ref{sec_prelims_block_graphs}, \ref{sec_block_graph}, and \ref{sec_combinatorial_char_reg_for_block_graphs}, we take up the question of understanding $\reg(J_{G})$ in the case when $G$ is a block graph. Previously, various authors have considered the case when $G$ is a tree (a special case of a block graph); see, for instance, Jayanthan et al.\ \cite{jayanthan2019regularityBlock} and the references therein. In \cite{herzog2018extremal}, Herzog and Rinaldo studied the extremal Betti numbers of $J_{G}$ when $G$ is a block graph and obtained a combinatorial characterization for $\reg(J_{G})$ for a subclass of all block graphs. The question of providing a combinatorial description of $\reg(J_{G})$ when $G$ is a block graph has remained open and was singled out by Herzog and Rinaldo as an important open question in the theory of binomial edge ideals \cite{herzog2018extremal}.

In Section \ref{sec_combinatorial_char_reg_for_block_graphs}, we prove the main result of this paper:
\begin{theoremx}
For a block graph $G$,
\begin{enumerate}
    \item $\nu(G)$ admits a combinatorial characterization solely in terms of vertex-disjoint induced paths of $G$ that do not admit certain induced subgraphs (Theorem \ref{thm_doip_iff_no_forbidden_path}),
    \item $\nu(G) = \reg(J_{G})-1$ (Theorem \ref{thm_reg=nu_block_graph}).
\end{enumerate}
\end{theoremx}
Theorem \ref{thm_doip_iff_no_forbidden_path} does not depend on \textit{any} choice of labeling of the induced paths nor on a choice of labeling of the block graph $G$. We prove Theorem \ref{thm_reg=nu_block_graph} by adapting the theory of Malayeri--Saeedi Madani--Kiani developed in \cite{malayeri2021proof}, where they provide a method to check whether a function is an upper bound for $\reg(J_{G})-1$ for every graph $G$. This answers the question of Herzog and Rinaldo. Moreover, this work shows that $\nu(G)$ gives a uniform computation for $\reg(J_{G})-1$ across many of the families of graphs considered thus far in the literature.

\section{Background} \label{sec_background}

\subsection{Graphs}

A (multi)graph $G$ is a pair $(V,E)$ where $V$ is a set and the elements are called vertices and $E$ is a multiset of pairs of vertices $\{a,b\}$, where we possibly allow repetition of the vertices appearing in an edge. When we wish to emphasize the vertex set (respectively edge set) of $G$, we will write $V(G)$ (respectively $E(G)$). An element appearing in $E$ multiple times is called a multi-edge, and an element of $E$ of the form $\{v,v\}$ for some $v \in V$ is called a loop. A graph having no loops nor multi-edges is called a simple graph, whereas a graph potentially having loops or multi-edges is called a multigraph. In this paper, when we write `graph' without the adjective `simple' or `multigraph,' we will implicitly mean a simple graph. When we wish to consider a multigraph, we will explicitly state that the graph is a multigraph. By a labeling of a set of vertices $W \subset V$, we mean a choice of an injective map $\phi : W \ra S$ where $S$ is a set of labels. When the vertices of $G$ have been labeled by a set possessing a total order, we will utilize the notation $v < w$ for vertices $v$ and $w$ to mean that $\phi(v) < \phi(w)$, where $\phi$ is the choice of labeling. By $[n]$, we denote the set of integers from $1$ to $n$ inclusive. By a graph $G$ on $[n]$, we mean that $\card{V(G)} = n$, and there is a labeling of the vertices of $G$ by $[n]$.

For a vertex $v \in V(G)$, we define the \textbf{neighbors} of $v$ in $G$ as the set: 
\begin{align*}
    N_{G}(v) := \{w \mid \{v,w\} \in  E(G)\},
\end{align*}
and we define the \textbf{degree} of $v$ in $G$ by
\begin{align*}
    \deg_{G}(v) := \card{N_{G}(v)}.
\end{align*}

We recall the following graph-theoretic constructions. For further information on the terminology introduced here, we refer the reader to \cite{westIntroToGraphTheory}. For a graph $G$ and $W \subset V(G)$, we define the \textbf{induced subgraph of $G$ on $W$}, which we denote by $\Ind_{G}(W)$ or by $G[W]$, as follows:
\begin{align*}
    V(\Ind_{G}(W)) &:= W \\
    E(\Ind_{G}(W)) &:= \{ \{a,b\} \mid \{a,b\} \in E(G), a\in W, b \in W \}.
\end{align*}
Given a subgraph $H$ of $G$, we will say that $H$ is an \textbf{induced subgraph} of $G$ if $\Ind_{G}(H) = H$. For a graph $G$ and $W \subset V(G)$, we define the graph $G \setminus W$ as follows:
\begin{align*}
    V(G\setminus W) &:= V(G) \setminus W \\
    E(G\setminus W) &:= \{ \{a,b\} \in E(G) \mid a\notin W, b \notin W\}.
\end{align*}
For a connected graph $G$ and $v \in V(G)$, we say that $v$ is a \textbf{cut vertex} of $G$ if $G \setminus \{v\}$ has strictly more connected components than $G$. 

We recall the well-known result that being an induced subgraph is transitive.
\blem
\label{lem_transitivity_of_being_induced}
Let $K$ be an induced subgraph of $H$, and $H$ be an induced subgraph of $G$. Then, $K$ is an induced subgraph of $G$.
\elem

For a graph $G$, we define a \textbf{path} of $G$ to be either: (i) a singleton vertex $v$ having no edges, or (ii) a sequence of vertices and edges $v_{1},e_{1},v_{2},e_{2},\ldots,v_{n-1},e_{n-1},v_{n}$ for some $n \geq 1$ satisfying:
\begin{enumerate}
    \item $v_{i} \in V(G)$,
    \item $e_{i} \in E(G)$,
    \item $v_{i} \neq v_{j}$ for all $i \neq j$, and 
    \item $e_{j} = \{v_{j},v_{j+1}\}$ for all $1 \leq j \leq n-1$.
\end{enumerate}
In the above setup, we will denote the above path via the notation $[v_{1},\ldots,v_{n}]$. For a path $P$ of $G$ and $v \in V(P)$, we say that $v$ is a \textbf{terminal vertex} of $P$ if $\deg_{P}(v) = 1$. We will say that a vertex $v \in V(P)$ is an \textbf{internal vertex} if $\deg_{P}(v) = 2$. We will denote by $\d P$ the set of terminal vertices of $P$, and we will denote by $P^{\circ}$ the set of internal vertices of $P$.

\bdefn
\label{defn_induced_edge}
Let $P_{i}$ be an induced path of $G$ for $1 \leq i \leq \ell$, and suppose that $V(P_{i}) \cap V(P_{j}) = \varnothing$ for all $1 \leq i < j \leq \ell$. An edge $e = \{a,b\} \in E(G)$ will be called an \textbf{induced edge} with respect to $P_{1},\ldots,P_{\ell}$ if $a\in P_{i}$ and $b \in P_{j}$ for some $1 \leq i \neq j \leq \ell$. We will denote by $P_{\Ind}$ the induced subgraph of $G$ on the vertices $\bigcup_{i=1}^{\ell} V(P_{i})$. A vertex $v \in V(P_{\ind})$ will be called an \textbf{internal vertex} (respectively, \textbf{terminal vertex}) of $P_{\ind}$ whenever $v$ in an internal vertex (respectively, terminal vertex) of $P_{i}$ for some $1 \leq i \leq \ell$. 
\edefn

We recall the definition of a directed graph and a key lemma about directed acyclic graphs.

\bdefn  
A \textbf{directed graph} $G$ consists of a set of vertices $V$ and a set of directed edges (or arcs) $A$. We denote a \textbf{directed edge} of a graph $G$ from the vertex $v$ to the vertex $w$ (with $v \neq w$) by $(v,w)$. Given a directed edge $(v,w)$, we say that $v$ and $w$ are the \textbf{initial vertex} and \textbf{terminal vertex}, respectively. When we wish to emphasize the set of arcs associated to $G$, we will utilize the notation $A(G)$. 

We will write \textbf{directed multigraph} to indicate that repetitions of directed edges are allowed and that directed edges having the same initial and terminal vertex are allowed. We will explicitly state when an object is a directed multigraph; otherwise, by directed graph, we always assume that there are no loops or repetitions of directed edges.

We say that a directed (multi)graph $G$ is \textbf{directed acyclic} if $G$ does not contain any directed cycle or any loop. A \textbf{topological ordering} or \textbf{topological sorting} of a directed (multi)graph $G$ is an integer labeling of the vertices of $G$ such that whenever $(i,j)$ is a directed edge of $G$, then $j < i$.
\edefn 

\blem[{\cite[Theorem 5.13, p.118]{thulasiraman2011graphs}}]
\label{lem_directed_acyclic_equiv_topological_sorting}
Let $G$ be a directed graph. $G$ is directed acyclic if and only if $G$ admits a topological sorting.
\elem


\subsection{Binomial Edge Ideals}
The main object of study in this paper are binomial edge ideals, introduced by Herzog et al.\ \cite{herzog2010binomial} and Ohtani \cite{ohtani2011graphs}, which associate to any simple graph a binomial ideal as follows. For a survey of binomial edge ideals, the reader is referred to \cite{saeedi2016binomial} or \cite{mukherjee2022homological}.

\bdefn[\cite{herzog2010binomial}]
\label{defn_bin_edge_ideal}
Let $G = (V,E)$ be a finite simple graph with vertex set $V$ labeled by $\{1,\ldots,n\}$ and edge set $E$. Fix a field $\K$. Consider the polynomial ring $S := \K[x_{1},\ldots,x_{n},y_{1},\ldots,y_{n}]$, and for each edge $\{i,j\} \in E$ with $i < j$ define $f_{ij} := x_{i}y_{j} - x_{j}y_{i} \in S$. Define the binomial edge ideal of $G$, denoted $J_{G}$, to be the ideal
\begin{equation}
\label{eqn_binomial_edge_ideal}
J_{G} := ( \{ f_{ij} \mid \{i,j\} \in E \} ).
\end{equation}
\edefn

In \cite{herzog2010binomial}, the authors provided a combinatorial description for a Gr\"obner basis of $J_{G}$ with respect to the lexicographic term order on $S$ induced by $x_{1} > x_{2} > \cdots > x_{n}> y_{1} > y_{2} > \cdots > y_{n}$. Throughout this paper, we will only consider this term order on $S$. We recall their result below.

\bdefn[{\cite[p.6]{herzog2010binomial}}]
\label{defn_adm_path}
Let $G$ be a simple graph on $[n]$, and let $i$ and $j$ be two vertices of $G$ with $i < j$. A path on the vertices $i_{0},i_{1},\ldots,i_{r}$ of $G$ with $i = i_{0}$ and $i_{r} = j$ is called \textbf{admissible} if:
\begin{enumerate}
\item $i_{k} \neq i_{l}$, for all $1 \leq k \neq l \leq r$, \label{item_adm_1}
\item for each $k = 1,\ldots,r-1$ one has either $i_{k} < i$ or $i_{k} > j$, \label{item_adm_2}
\item for any proper subset $\{j_{1},\ldots,j_{s}\}$ of $\{i_{1},\ldots,i_{r-1}\}$ the sequence $i,j_{1},\ldots,j_{s},j$ is not a path. \label{item_adm_3}
\end{enumerate}
Given such an admissible path, we define the monomial
\begin{align*}
u_{\pi} = \left( \prod_{i_{k}>j} x_{i_{k}} \right) \left( \prod_{i_{l} < i} y_{i_{l}} \right),
\end{align*}
and we denote by $m_{\pi}$ the monomial $x_{i}y_{j}u_{\pi}$.
\edefn

\brem
\cref{item_adm_1} and \cref{item_adm_3} of Definition \ref{defn_adm_path} establish that $\pi$ is an induced path of $G$.
\erem

\bthm[{\cite[Theorem 2.1]{herzog2010binomial}}]
\label{thm_gro_basis}
Let $G$ be a simple graph on $[n]$. Let $<$ be the lexicographic order on $S = \K[x_{1},\ldots,x_{n},y_{1},\ldots,y_{n}]$ induced by $x_{1} > x_{2} > \cdots > x_{n}> y_{1} > y_{2} > \cdots > y_{n}$. Then, the set of binomials
\begin{align*}
\mathscr{G} := \bigcup_{i<j} \, \{u_{\pi} f_{ij} \mid \pi \text{ is an addmissible path from } i \text{ to } j\}
\end{align*}
is a reduced Gr\"{o}bner basis of $J_{G}$. 
\ethm
Consequently, $J_{G}$ is a radical ideal {\cite[Corollary 2.2]{herzog2010binomial}}.

\subsection{Castelnuovo--Mumford Regularity}

We recall the definition of Castelnuovo--Mumford regularity of a finitely generated graded $R$-module, where $R$ is a polynomial ring. Let $R := \K[z_{1},\ldots,z_{m}]$ be standard graded. Given a finitely generated graded $R$-module $M$, let 
\begin{align*}
F_{\bullet}: 0 \ra F_{n} \ra F_{n-1} \ra \cdots \ra F_{1} \ra F_{0} \ra 0
\end{align*}
be a minimal graded free $R$-resolution of $M$ where $F_{i} = \bigoplus_{j\in \Z} R(-j)^{b_{ij}}$. The $b_{ij}$ are the \textbf{graded Betti numbers} of $M$, non-negative integers, and for each $i$, only finitely many of the $b_{ij}$ are non-zero. The \textbf{Castelnuovo--Mumford regularity} of $M$ is defined as follows:
\begin{equation}
\label{defn_regularity}
\reg(M) := \max\{ j-i \mid b_{ij} \neq 0\}.
\end{equation}
The reader is referred to \cite{peeva2010graded} or \cite{bruns2021castelnuovo} for further information regarding Castelnuovo-Mumford regularity.

The next result of Conca and Varbaro \cite{conca2020square} shows that under the assumption that a homogeneous ideal has a squarefree initial ideal (with respect to some term order), then the extremal Betti numbers of the ideal and of its initial ideal coincide; in particular, their regularities coincide.

\bthm[{\cite[Corollary 2.7]{conca2020square}}]
\label{thm_conca_varbaro}
Let $I \subset R := \K[z_{1},\ldots,z_{m}]$ be a homogeneous ideal such that $\init(I)$ is square-free with respect to some term order (not necessarily lexicographic order). Then, the extremal Betti numbers of $R/I$ and those of $R/\init(I)$ coincide (positions and values). In particular, $\reg(R/I) = \reg(R/\init(I))$.
\ethm

\section{The Invariant $\nu(G)$} \label{sec_the_invariant_nu}

\subsection{Definition and Motivation}

It is a result of Matsuda and Murai that:
\bthm[{\cite[Corollary 2.2]{matsuda2013regularity}}]
\label{thm_mm_induced_subgraph_bound_reg}
If $H$ is an induced subgraph of $G$, then $\reg(S/J_{H}) \leq \reg(S/J_{G})$.
\ethm
Theorem \ref{thm_mm_induced_subgraph_bound_reg} is most often used in the form of the following corollary.
\bcor[{\cite[Corollary 2.3]{matsuda2013regularity}}]
\label{cor_mm_induced_path_leq_reg}
Let $G$ be a graph, then 
\begin{align*}
    \ell(G)\leq \reg(S/J_{G})
\end{align*}
where $\ell(G)$ is the length of a longest induced path within $G$.
\ecor

However, Theorem \ref{thm_mm_induced_subgraph_bound_reg} also implies the slightly stronger result that if $P_{1}, \ldots,P_{\ell}$ are vertex-disjoint induced paths of $G$ having no induced edges (Definition \ref{defn_induced_edge}), then $\sum_{i=1}^{\ell} \card{E(P_{i})} \leq~\reg(S/J_{G})$. It is perhaps natural to ask:
\bques
\label{ques_allowing_some_induced_edges}
What induced edges can we allow between vertex-disjoint induced paths $P_{1},\ldots,P_{\ell}$ of $G$ while retaining the lower bound 
\begin{equation}
    \label{eqn:lb}
    \sum_{i=1}^{\ell} \card{E(P_{i})} \leq \reg(S/J_{G})?
\end{equation}
\eques
For example, it can be checked with Macaulay 2 \cite{M2} that the graph $G$ in Figure \ref{fig_doip_graph}, consisting of the induced paths $[4,6,5]$ and $[1,3,2]$ together with the induced edge $\{3,4\}$, satisfies $\reg(S/J_{G}) = 4$. 

\def\putLeft{.5*.9}
\def\edgeDist{1.75*.9}
\begin{figure}
\begin{center}
\begin{tikzpicture}
\filldraw[black] (0*\edgeDist,0*\edgeDist) circle (2pt) node at (-1*\putLeft,0*\edgeDist) {$5$};
\filldraw[black] (0*\edgeDist,1*\edgeDist) circle (2pt) node at (-1*\putLeft,1*\edgeDist) {$6$};
\filldraw[black] (0*\edgeDist,2*\edgeDist) circle (2pt) node at (-1*\putLeft,2*\edgeDist) {$4$};

\filldraw[black] (1*\edgeDist,0*\edgeDist) circle (2pt) node at (1*\edgeDist+\putLeft,0*\edgeDist) {$2$};
\filldraw[black] (1*\edgeDist,1*\edgeDist) circle (2pt) node at (1*\edgeDist+\putLeft,1*\edgeDist) {$3$};
\filldraw[black] (1*\edgeDist,2*\edgeDist) circle (2pt) node at (1*\edgeDist+\putLeft,2*\edgeDist) {$1$};

\draw (0*\edgeDist,0*\edgeDist) to (0*\edgeDist,1*\edgeDist); 
\draw (0*\edgeDist,1*\edgeDist) to (0*\edgeDist,2*\edgeDist); 
\draw (0*\edgeDist,2*\edgeDist) to (1*\edgeDist,1*\edgeDist); 
\draw (1*\edgeDist,0*\edgeDist) to (1*\edgeDist,1*\edgeDist); 
\draw (1*\edgeDist,1*\edgeDist) to (1*\edgeDist,2*\edgeDist); 

\end{tikzpicture}
\end{center}
\caption{DOIP Paths}
\label{fig_doip_graph}
\end{figure}

\subsection{Directed Oriented Induced Paths}

We introduce the following definition, which provides a sufficient condition for vertex-disjoint induced paths to satisfy equation \eqref{eqn:lb}.

\bdefn
\label{defn_oriented_path}
Let $P$ be an induced path of a graph $G$. We say that $P$ is an \textbf{oriented induced path} if there exists a surjection $\phi : \{1,2\} \ra \d P$. We will refer to $\phi$ as an \textbf{orientation}.
\edefn 

\brem
When an oriented path $P$ has exactly one vertex, $\phi_{P}(1) = \phi_{P}(2)$. Otherwise, $\phi_{P}$ is a bijection. In the latter case, we think of $\phi$ as specifying a start and an end vertex for $P_{i}$. This distinction between the terminal vertices of $P_{i}$ will prove necessary due to the asymmetry between the terminal vertices of admissible paths in Definition \ref{defn_adm_path}.
\erem 

\bdefn
\label{defn_doip}
Let $\ul{P} := P_{1},\ldots,P_{\ell}$ be vertex-disjoint induced paths of $G$. Let $P_{\ind}$ denote the induced subgraph of $G$ on $\bigcup_{i=1}^{\ell} V(P_{i})$. For a choice of
\begin{enumerate}
    \item $\sigma$ a permutation on the set $\{1,\ldots,\ell\}$, and
    \item orientations $\phi_{i} : \{1,2\} \ra \d P_{i}$ for $1 \leq i \leq \ell$,
\end{enumerate}
we say that $(\ul{P},\sigma,\ul{\phi_{i}})$ are \textbf{directed oriented induced paths (DOIP)} if whenever $\sigma(i) \leq \sigma(j)$ and $Q$ is an induced path of $P_{\ind}$ having terminal vertices $\phi_{\sigma(i)}(1)$ and $\phi_{\sigma(j)}(2)$, then $Q$ contains $P_{k}$ as subgraph for some $1 \leq k \leq \ell$. We will say that $\ul{P}$ is DOIP if there exists a choice of $\sigma$ and orientations $\ul{\phi_{i}}$ such that $(\ul{P},\sigma,\ul{\phi_{i}})$ is DOIP.
\edefn 

\bex
\label{ex_doip}
In Figure \ref{fig_doip_graph}, we define the paths $P_{1} = [1,3,2]$ and $P_{2} = [4,6,5]$. Then, $P_{1}$ and $P_{2}$ are DOIP. Indeed, we let $\sigma = \id_{\{1,2\}}$, and $\phi_{1}(1) = 1$, $\phi_{1}(2) = 2$, $\phi_{2}(1) = 4$, and $\phi_{2}(2) = 5$. It is now clear that any induced path of $P_{\ind}$ from vertex $1$ to either vertex $2$ or vertex $4$ contains either $P_{1}$ or $P_{2}$. Likewise, for any induced path from vertex $4$ to vertex $5$.
\eex 

\brem
If in Example \ref{ex_doip}, we were to change $\sigma$ from the identity permutation to the transposition $(2 \, 1)$ while keeping $P_{1}$, $P_{2}$, $\phi_{1}$, and $\phi_{2}$ as in Example \ref{ex_doip}, then $[4,3,2]$ is an induced path of $P_{\ind}$ from $\phi_{\sigma(1)}(1)$ to $\phi_{\sigma(2)}(2)$ which does not contain $P_{1}$ or $P_{2}$. 

If in Example \ref{ex_doip}, we were to keep $P_{1}$, $P_{2}$, and $\phi_{1}$ unchanged, and we were to change $\sigma$ from the identity permutation to the transposition $(2 \, 1)$ and $\phi_{2}$ to $\phi_{2}(1) = 5$ and $\phi_{2}(2) = 4$, then $P_{1}$ and $P_{2}$ are DOIP with respect to these choices. These examples demonstrate that the property of $(\ul{P},\sigma,\ul{\phi_{i}})$ being DOIP may depend on $\sigma$ and the $\phi_{i}$.
\erem 

In the next example, we demonstrate vertex-disjoint induced paths which are not DOIP.

\bex
We consider Figure \ref{fig_non_doip}. Consider the induced paths $P_{1} = [1,3,2]$ and $P_{2} = [4,6,5]$ in each of the three graphs depicted in this figure. Then, $P_{1}$ and $P_{2}$ are not DOIP. For each of these graphs, any pair of terminal vertices from $P_{1}$ and $P_{2}$ can be connected by an induced path not containing $P_{1}$ or $P_{2}$. The existence of such paths is an obstruction to the paths $P_{1}$ and $P_{2}$ being DOIP.

Furthermore, for the center and rightmost graphs, there is an induced path connecting the terminal vertices of $P_{1}$, which does not contain $P_{1}$ nor $P_{2}$. The existence of such a path is also an obstruction to the DOIP property.

We observe that we can find subpaths in these graphs which are DOIP; the paths $P_{1} = [4,6]$ and $P_{2} = [1,3,2]$ are DOIP for all of these graphs.
\eex 

\def\putLeft{.5*.9}
\def\edgeDist{1.75*.9}
\begin{figure}
\begin{subfigure}{.3\linewidth}
\begin{tikzpicture}
\filldraw[black] (0*\edgeDist,0*\edgeDist) circle (2pt) node at (-1*\putLeft,0*\edgeDist) {$5$};
\filldraw[black] (0*\edgeDist,1*\edgeDist) circle (2pt) node at (-1*\putLeft,1*\edgeDist) {$6$};
\filldraw[black] (0*\edgeDist,2*\edgeDist) circle (2pt) node at (-1*\putLeft,2*\edgeDist) {$4$};

\filldraw[black] (1*\edgeDist,0*\edgeDist) circle (2pt) node at (1*\edgeDist+\putLeft,0*\edgeDist) {$2$};
\filldraw[black] (1*\edgeDist,1*\edgeDist) circle (2pt) node at (1*\edgeDist+\putLeft,1*\edgeDist) {$3$};
\filldraw[black] (1*\edgeDist,2*\edgeDist) circle (2pt) node at (1*\edgeDist+\putLeft,2*\edgeDist) {$1$};

\draw (0*\edgeDist,0*\edgeDist) to (0*\edgeDist,1*\edgeDist); 
\draw (0*\edgeDist,1*\edgeDist) to (0*\edgeDist,2*\edgeDist); 
\draw (0*\edgeDist,1*\edgeDist) to (1*\edgeDist,1*\edgeDist); 
\draw (1*\edgeDist,0*\edgeDist) to (1*\edgeDist,1*\edgeDist); 
\draw (1*\edgeDist,1*\edgeDist) to (1*\edgeDist,2*\edgeDist); 

\end{tikzpicture}
\end{subfigure}
%
\begin{subfigure}{.3\linewidth}
\begin{tikzpicture}
\filldraw[black] (0*\edgeDist,0*\edgeDist) circle (2pt) node at (-1*\putLeft,0*\edgeDist) {$5$};
\filldraw[black] (0*\edgeDist,1*\edgeDist) circle (2pt) node at (-1*\putLeft,1*\edgeDist) {$6$};
\filldraw[black] (0*\edgeDist,2*\edgeDist) circle (2pt) node at (-1*\putLeft,2*\edgeDist) {$4$};

\filldraw[black] (1*\edgeDist,0*\edgeDist) circle (2pt) node at (1*\edgeDist+\putLeft,0*\edgeDist) {$2$};
\filldraw[black] (1*\edgeDist,1*\edgeDist) circle (2pt) node at (1*\edgeDist+\putLeft,1*\edgeDist) {$3$};
\filldraw[black] (1*\edgeDist,2*\edgeDist) circle (2pt) node at (1*\edgeDist+\putLeft,2*\edgeDist) {$1$};

\draw (0*\edgeDist,0*\edgeDist) to (0*\edgeDist,1*\edgeDist); 
\draw (0*\edgeDist,1*\edgeDist) to (0*\edgeDist,2*\edgeDist); 
\draw (0*\edgeDist,2*\edgeDist) to (1*\edgeDist,0*\edgeDist); 
\draw (0*\edgeDist,0*\edgeDist) to (1*\edgeDist,2*\edgeDist); 
\draw (0*\edgeDist,2*\edgeDist) to (1*\edgeDist,2*\edgeDist); 
\draw (0*\edgeDist,0*\edgeDist) to (1*\edgeDist,0*\edgeDist); 
\draw (1*\edgeDist,0*\edgeDist) to (1*\edgeDist,1*\edgeDist); 
\draw (1*\edgeDist,1*\edgeDist) to (1*\edgeDist,2*\edgeDist); 

\end{tikzpicture}
\end{subfigure}
%
\begin{subfigure}{.3\linewidth}
\begin{tikzpicture}
\filldraw[black] (0*\edgeDist,0*\edgeDist) circle (2pt) node at (-1*\putLeft,0*\edgeDist) {$5$};
\filldraw[black] (0*\edgeDist,1*\edgeDist) circle (2pt) node at (-1*\putLeft,1*\edgeDist) {$6$};
\filldraw[black] (0*\edgeDist,2*\edgeDist) circle (2pt) node at (-1*\putLeft,2*\edgeDist) {$4$};

\filldraw[black] (1*\edgeDist,0*\edgeDist) circle (2pt) node at (1*\edgeDist+\putLeft,0*\edgeDist) {$2$};
\filldraw[black] (1*\edgeDist,1*\edgeDist) circle (2pt) node at (1*\edgeDist+\putLeft,1*\edgeDist) {$3$};
\filldraw[black] (1*\edgeDist,2*\edgeDist) circle (2pt) node at (1*\edgeDist+\putLeft,2*\edgeDist) {$1$};

\draw (0*\edgeDist,0*\edgeDist) to (0*\edgeDist,1*\edgeDist); 
\draw (0*\edgeDist,1*\edgeDist) to (0*\edgeDist,2*\edgeDist); 
\draw (0*\edgeDist,2*\edgeDist) to (1*\edgeDist,1*\edgeDist); 
\draw (0*\edgeDist,0*\edgeDist) to (1*\edgeDist,1*\edgeDist); 
\draw (1*\edgeDist,0*\edgeDist) to (1*\edgeDist,1*\edgeDist); 
\draw (1*\edgeDist,1*\edgeDist) to (1*\edgeDist,2*\edgeDist); 

\end{tikzpicture}
\end{subfigure}
\caption{Non-DOIP paths}
\label{fig_non_doip}
\end{figure}

The notion of $P_{1},\ldots,P_{\ell}$ being DOIP captures the idea that any induced path $Q$ of $P_{\ind}$ with $\d Q \subset \bigcup_{i=1}^{\ell} \d P_{i}$, such that $Q$ does not contain some $P_{k}$, \textit{travels} from top to bottom and from left to right. We make this idea precise using the notion of directed acyclic graphs.

\bdefn
\label{defn_turning_admissible_blocking_paths_into_directed_graph}
Let $\ul{P} := P_{1},\ldots,P_{\ell}$ be vertex-disjoint induced paths of $G$, let $P_{\ind}$ be the induced subgraph of $G$ on $\bigcup_{i=1}^{\ell} V(P_{i})$, and let $\phi_{i} : \{1,2\} \ra \d P_{i}$ be orientations for $1 \leq i \leq \ell$. Then, we define $K_{P_{\ind}}$ to be the directed multigraph with vertex set $[\ell]$, and with a multiarc $(i,j)$ for each induced path $Q$ from $\phi_{i}(1)$ to $\phi_{j}(2)$ whenever $1 \leq i,j \leq \ell$ and $Q$ does not contain $P_{k}$ for every $1 \leq k \leq \ell$. 
\edefn 

\brem 
Up to isomorphism of multigraphs, $K_{P_{\ind}}$ does not depend on the choice of labeling $\sigma$ of the paths $P_{1},\ldots,P_{\ell}$. However, as Example \ref{ex_dep_K_P_on_orientations} illustrates, $K_{P_{\ind}}$ does depend on the choice of orientations $\phi_{i}$.
\erem 

\bex
\label{ex_dep_K_P_on_orientations}
In Figure \ref{fig_K_P}, let $P_{1} = [1,2]$ and $P_{2} = [3,5,4]$. Let $\phi_{1}(1) = 1$, $\phi_{1}(2) = 2$, $\phi_{2}(1) = 3$, and $\phi_{2}(2) = 4$. Then, $K_{P_{\ind}}$ is the directed multigraph on the vertex set $\{1,2\}$ with multiarcs: $(2,1)$ corresponding to the induced path $[3,5,2]$.

Now, let us suppose that $\phi_{1}(1) = 2$, $\phi_{1}(2) = 1$, $\phi_{2}(1) = 3$, and $\phi_{2}(2) = 4$. Then, $K_{P_{\ind}}$ is the directed multigraph on the vertex set $\{1,2\}$ with multiarcs: $(2,1)$ corresponding to the induced path $[3,1]$, and $(1,2)$ corresponding to the induced path $[2,5,4]$.
\eex 

\def\putLeft{.5*.9}
\def\edgeDist{1.75*.9}
\begin{figure}
\begin{center}
\begin{tikzpicture}
\filldraw[black] (0*\edgeDist,0*\edgeDist) circle (2pt) node at (-1*\putLeft,0*\edgeDist) {$4$};
\filldraw[black] (0*\edgeDist,1*\edgeDist) circle (2pt) node at (-1*\putLeft,1*\edgeDist) {$5$};
\filldraw[black] (0*\edgeDist,2*\edgeDist) circle (2pt) node at (-1*\putLeft,2*\edgeDist) {$3$};

\filldraw[black] (1*\edgeDist,1*\edgeDist) circle (2pt) node at (1*\edgeDist+\putLeft,1*\edgeDist) {$2$};
\filldraw[black] (1*\edgeDist,2*\edgeDist) circle (2pt) node at (1*\edgeDist+\putLeft,2*\edgeDist) {$1$};

\draw (0*\edgeDist,0*\edgeDist) to (0*\edgeDist,1*\edgeDist); 
\draw (0*\edgeDist,1*\edgeDist) to (0*\edgeDist,2*\edgeDist); 
\draw (0*\edgeDist,1*\edgeDist) to (1*\edgeDist,1*\edgeDist); 
\draw (1*\edgeDist,1*\edgeDist) to (1*\edgeDist,2*\edgeDist); 
\draw (0*\edgeDist,2*\edgeDist) to (1*\edgeDist,2*\edgeDist); 

\end{tikzpicture}
\end{center}
\caption{Dependence of $K_{P_{\ind}}$ on Orientations}
\label{fig_K_P}
\end{figure}

\bthm
\label{thm_doip_equiv_K_P_acyclic}
Let $\ul{P} := P_{1},\ldots,P_{\ell}$ be vertex-disjoint induced paths of $G$. Then, $\ul{P}$ is DOIP if and only if there exist orientations $\phi_{i} : \{1,2\} \ra \d P_{i}$ for $1 \leq i \leq \ell$ such that the multigraph $K_{P_{\ind}}$ corresponding to these orientations is directed acyclic.
\ethm 

\bproof
$(\implies)$ Immediate consequence of Definition \ref{defn_doip}.

$(\impliedby)$ By Lemma \ref{lem_directed_acyclic_equiv_topological_sorting}, there exists an ordering of the vertices of $K_{P_{\ind}}$ under which $K_{P_{\ind}}$ admits a topological sorting. Let $\sigma$ be the bijection which realizes this ordering of the vertices of $K_{P_{\ind}}$.
\eproof 

\bcor 
\label{cor_induced_path_is_doip}
An induced path $P_{1}$ of $G$ is DOIP.
\ecor 

\bproof
Follows immediately from Theorem \ref{thm_doip_equiv_K_P_acyclic}, since $K_{P_{\ind}}$ is a singleton vertex possessing no multiarcs.
\eproof 

\subsection{DOIP and Regularity} In this subsection, we establish that paths which are DOIP satisfy equation \eqref{eqn:lb}.

\bdefn 
\label{defn_supp_monl_and_ideal}
Let $m \in S$ be a monomial. We define the \textbf{support} of $m$ to be the subset of variables of $S$ which divide $m$. We denote by $\Supp(m)$ the support of $m$. For a monomial ideal $I \subset S$ and $W \subset S$ any set of monomials, we define
\begin{align*}
    I_{W} := (\{m \in I \mid m \text{ is a monomial and } \Supp(m) \subset W\}).
\end{align*}
\edefn 

\blem
\label{lem_doip_implies_initial_ideal_good}
Let $\ul{P} := P_{1},\ldots,P_{\ell}$ be vertex-disjoint induced paths of a graph $G$ which are DOIP. Then, there exists a labeling of the vertices of $G$ such that $P_{i}$ with respect to this labeling is an admissible path in the sense of Definition \ref{defn_adm_path} for $1 \leq i \leq \ell$, and that
\begin{equation}
    \label{eqn_doip_implies_initial_ideal_good_1}
    \left( \init (J_{G}) \right)_{W} = ( m_{1},\ldots,m_{\ell} ).
\end{equation}
Where we denote the monomial associated to $P_{i}$ in Definition \ref{defn_adm_path} by $m_{i}$, and we define $W := \bigcup_{i=1}^{\ell} \Supp(m_{i})$.
\elem 

\bproof
We may suppose that the paths $P_{1},\ldots,P_{\ell}$ are labeled such that $(\ul{P},\id_{[\ell]},\ul{\phi_{i}})$ is DOIP. For $1 \leq i \leq \ell$, label vertex $\phi_{i}(1)$ by the integer $2i-1$, and label vertex $\phi_{i}(2)$ by the integer $2i$. Label the remaining vertices of $G$ by distinct consecutive integers larger than $2\ell$. With respect to this labeling, the $P_{i}$ are admissible for each $1 \leq i \leq \ell$. Moreover, by Theorem \ref{thm_gro_basis}, we have that
\begin{align*}
    y_{j} \in W := \bigcup_{i=1}^{\ell} \Supp(m_{i})
\end{align*}
if and only if $j = 2i$ for $1 \leq i \leq \ell$.

We next establish equation \eqref{eqn_doip_implies_initial_ideal_good_1}. As the reverse inclusion is clear, it suffices to prove the forward inclusion. Let $m \in \init (J_{G})$ be a monomial such that $\Supp(m) \subset W$. Then, there exist monomials $u \in S$ and $m^{'} \in \init(J_{G})$, a minimal generator, such that $m = u\cdot m^{'}$. As Theorem \ref{thm_gro_basis} gives a reduced Gr\"obner basis of $J_{G}$, we have that $m^{'} = m_{Q}$ for some admissible path $Q$ of $G$. As $V(Q) \subset W$, it follows, in particular, that $Q$ is an induced path of $P_{\ind}$. In order for $\Supp(m_{Q}) \subset W$, it is necessary that one of the terminal vertices of $Q$ is $\phi_{i}(2)$ for some $1 \leq i \leq \ell$. In order for $Q$ to be admissible, it is necessary that the other terminal vertex of $Q$ is $\phi_{j}(1)$ for some $1 \leq j \leq i$ (because all the other vertices of $G$ are labeled by integers larger than $2i$). Now, as $\ul{P}$ is DOIP, Definition \ref{defn_doip} implies that $P_{k}$ is a subgraph of $Q$ for some $1 \leq k \leq \ell$. We observe that:
\begin{enumerate}
    \item $Q$ does not contain $P_{r}$ as a subgraph for $1 \leq r < j$. Otherwise, $y_{2r-1}\mid m_{Q}$, but $y_{2r-1} \notin W$ ($\because x_{2r-1} \in W$).
    \item $Q$ does not contain $P_{r}$ as a subgraph for $j < r < i$. Otherwise, $Q$ contains the vertex $2r-1$. But $2r-1$ is strictly between the terminal vertices of $Q$, which are $2j-1$ and $2i$, contradicting $Q$ being admissible.
    \item $Q$ does not contain $P_{r}$ as a subgraph for $i < r \leq \ell$. Otherwise, $x_{2r} \mid m_{Q}$, but $x_{2r} \notin W$ ($\because y_{2r} \in W$).
\end{enumerate}
It follows from these observations that $j = k = i$. In order for $Q$ to:
\begin{enumerate}
    \item contain $P_{k}$ as a subgraph, 
    \item be admissible, and
    \item have terminal vertices $2j-1$ and $2i$,
\end{enumerate}
it must be the case that $Q = P_{k}$. Consequently, $m \in ( m_{1},\ldots,m_{\ell} )$, which completes the proof.
\eproof

\bex 
We illustrate Lemma \ref{lem_doip_implies_initial_ideal_good} in the context of the graph in Figure \ref{fig_doip_graph}.

Let $P_{1} = [1,3,2]$, and $P_{2} = [4,6,5]$. We observe that $m_{1} = x_{1}x_{3}y_{2}$ and that $m_{2} = x_{4}x_{6}y_{5}$. Hence, $W = \{x_{1},x_{3},x_{4},x_{6},y_{2},y_{5}\}$. It can be checked that 
\begin{align*}
    \init(J_{G}) = ( x_5 y_6, x_4 y_6, x_4 x_6 y_5, x_3 y_4, x_2 y_3, x_1 y_3, x_1 x_3 y_2 ).
\end{align*}
We see that
\begin{align*}
    \init(J_{G})_{W} = ( x_4 x_6 y_5, x_1 x_3 y_2 ).
\end{align*}
\eex 

\bdefn 
\label{defn_nu}
Let $G$ be a graph, we define the invariant
\begin{align*}
    \nu(G) := \max \left \{ \sum_{i=1}^{\ell} \card{E(P_{i})} \mid \ul{P} \text{ is DOIP} \right \}.
\end{align*}
\edefn 

\brem
We observe that $\nu(G)$ can equivalently be expressed as:
\begin{align*}
    \nu(G) = \max_{\phi_{i} : 1 \leq i \leq \ell} \max_{\sigma} \max \left \{ \sum_{i=1}^{\ell} \card{E(P_{i})} \mid (\{P_{i}\}_{i=1}^{\ell},\sigma,\{\phi_{i}\}_{i=1}^{\ell}) \text{ is DOIP} \right \}.
\end{align*}
\erem 

\bthm 
\label{thm_nu_leq_reg}
Let $G$ be a graph, then
\begin{align*}
    \nu(G) \leq \reg(S/J_{G}).
\end{align*}
\ethm 

\bproof
It suffices to show that if $P_{1},\ldots,P_{\ell}$ is DOIP, then 
\begin{align*}
    \sum_{i=1}^{\ell} \card{E(P_{i})} \leq \reg(S/J_{G}).
\end{align*}
Label the vertices of $G$ as in Lemma \ref{lem_doip_implies_initial_ideal_good}. Since $\init(J_{G})$ is a squarefree monomial ideal, we have by Theorem \ref{thm_conca_varbaro} that
\begin{align*}
\reg(S/\init(J_{G})) = \reg(S/J_{G}).
\end{align*}
It is well known that
\begin{align*}
    \reg(S/\init(J_{G})_{W}) \leq \reg(S/\init(J_{G})).
\end{align*}
(See, for example, \cite{peeva2010graded}.) Lemma \ref{lem_doip_implies_initial_ideal_good} implies that $\init(J_{G})_{W} = (m_{1},\ldots,m_{\ell})$ is a complete intersection. Hence,
\begin{align*}
    \reg(S/(m_{1},\ldots,m_{\ell})) 
    &= \sum_{i=1}^{\ell} (\deg(m_{i}) -1) \\
    &= \sum_{i=1}^{\ell} \card{E(P_{i})},
\end{align*}
which completes the proof.
\eproof

We recover Matsuda and Murai's lower bound for Castelnuovo--Mumford regularity as a corollary.

\bcor[{\cite[Corollary 2.2]{matsuda2013regularity}}]
\label{cor_induced_path_leq_reg}
For a graph $G$, 
\begin{align*}
    \ell(G) \leq \nu(G) \leq \reg(S/J_{G}).
\end{align*}
\ecor 

\bproof
Follows immediately from Corollary \ref{cor_induced_path_is_doip} and Theorem \ref{thm_nu_leq_reg}.
\eproof 

\section{Equality of $\nu(G)$ and $\reg(S/J_{G})$} \label{sec_equality_nu_and_reg}

In this section, we recall two families of graphs for which a combinatorial description of the Castelnuovo--Mumford regularity of the binomial edge ideal is known. We show that for these two families, the previous computations for their Castelnuovo--Mumford regularity coincide with the invariant $\nu(G)$.

\subsection{Closed Graphs}

In \cite{herzog2010binomial}, closed graphs were introduced. A graph $G$ is said to be \textbf{closed} if there exists a labeling of the vertices such that the set $$\{ x_{i}y_{j} - x_{j}y_{i} \mid \{i,j\} \in E(G) \}$$ is a quadratic Gr\"obner basis of $J_{G}$ (\cite[Theorem 1.1]{herzog2010binomial}). Crupi and Rinaldo showed that closed graphs are interval graphs {\cite[Theorem 2.4]{crupi2014closed}}. In \cite{ene2015regularity}, Ene and Zarojanu computed the Castelnuovo--Mumford regularity for closed graphs.

\bthm[{\cite[Theorem 2.2]{ene2015regularity}}]
\label{thm_reg_closed_graph}
Let $G$ be a closed graph, then
\begin{align*}
    \ell(G) = \reg(S/J_{G}).
\end{align*}
\ethm 

\bcor
\label{cor_nu_equals_reg_closed_graphs}
Let $G$ be a closed graph, then
\begin{align*}
    \ell(G) = \nu(G) = \reg(S/J_{G}).
\end{align*}
\ecor

\bproof
Corollary \ref{cor_induced_path_is_doip} and Theorem \ref{thm_nu_leq_reg} imply that 
\begin{align*}
    \ell(G) \leq \nu(G) \leq \reg(R/J_{G}).
\end{align*}
Equality throughout now follows from Theorem \ref{thm_reg_closed_graph}.
\eproof

\subsection{Cohen--Macaulay Bipartite Graphs}

In \cite{bolognini2018binomial}, the authors study when the binomial edge ideal of a bipartite graph is unmixed, and they give a combinatorial characterization of when the binomial edge ideal of a bipartite graph is Cohen--Macaulay. Using this characterization of Cohen--Macaulay binomial edge ideals of bipartite graphs, Jayanthan and Kumar computed the regularity for this family of graphs {\cite[Theorem 4.7]{jayanthan2019regularityCM&BipartiteGraphs}}. We now recall this characterization of Cohen--Macaulay binomial edge ideals of bipartite graphs (we use the notation from both \cite{bolognini2018binomial} and \cite{jayanthan2019regularityCM&BipartiteGraphs}).

\bdefn[{\cite[p.2]{bolognini2018binomial}}]
\label{defn_operations_bipartite_graphs}
For every $m \geq 1$, let $F_{m}$ be the graph on the vertex set $[2m]$ and with edge set 
\begin{align*}
    E(F_{m}) := \{ (2i,2j-1)\} \mid i = 1,\ldots,m, j = i,\ldots,m\}.
\end{align*}

\ul{The operation $\ast$}: For $i = 1,2$, let $G_{i}$ be a graph having at least one vertex $f_{i}$ of degree one. We define $(G_{1},f_{1}) \ast (G_{2},f_{2})$ to be the graph obtained by identifying the vertices $f_{1}$ and $f_{2}$.

\ul{The operation $\circ$}: For $i = 1,2$, let $G_{i}$ be a graph with at least one vertex $f_{i}$ of degree one, and let $v_{i}$ be its neighbor, and we assume that $\deg_{G_{i}}(v_{i}) \geq 3$. We define $(G_{1},f_{1}) \circ (G_{2},f_{2})$ to be the graph obtained from $G_{1}$ and $G_{2}$ by deleting the vertices $f_{1}$, $f_{2}$ and identifying the vertices $v_{1}$ and $v_{2}$.
\edefn

Figure \ref{fig_Fm} depicts $F_{i}$ for $1 \leq i \leq 4$. Example \ref{ex_nu_equals_reg_fomula} and Figure \ref{fig_total_construction} includes a concrete description and illustration of these operations $\ast$ and $\circ$ for a particular bipartite graph $G$.

\def\putLeft{.5*.9}
\def\edgeDist{1.75*.9}

\begin{figure}
\begin{subfigure}{.1\textwidth} 
\centering
\begin{tikzpicture}
\filldraw[black] (0*\edgeDist,0*\edgeDist) circle (2pt) node at (0*\edgeDist,0*\edgeDist-\putLeft) {$2$};
\filldraw[black] (0*\edgeDist,1*\edgeDist) circle (2pt) node at (0*\edgeDist,1*\edgeDist+\putLeft) {$1$};
%
\draw (0*\edgeDist,0*\edgeDist) to (0*\edgeDist,1*\edgeDist); 
\end{tikzpicture}
\caption{$F_{1}$}
\end{subfigure}
\hfill
\begin{subfigure}{.15\textwidth}
\centering
\begin{tikzpicture}
\filldraw[black] (0*\edgeDist,0*\edgeDist) circle (2pt) node at (0*\edgeDist,0*\edgeDist-\putLeft) {$2$};
\filldraw[black] (0*\edgeDist,1*\edgeDist) circle (2pt) node at (0*\edgeDist,1*\edgeDist+\putLeft) {$1$};
\filldraw[black] (1*\edgeDist,0*\edgeDist) circle (2pt) node at (1*\edgeDist,0*\edgeDist-\putLeft) {$4$};
\filldraw[black] (1*\edgeDist,1*\edgeDist) circle (2pt) node at (1*\edgeDist,1*\edgeDist+\putLeft) {$3$};
%
\draw (0*\edgeDist,0*\edgeDist) to (0*\edgeDist,1*\edgeDist); 
\draw (0*\edgeDist,0*\edgeDist) to (1*\edgeDist,1*\edgeDist); 
\draw (1*\edgeDist,0*\edgeDist) to (1*\edgeDist,1*\edgeDist); 
\end{tikzpicture}
\caption{$F_{2}$}
\end{subfigure}
\hfill
\begin{subfigure}{.225\textwidth}
\centering
\begin{tikzpicture}
\filldraw[black] (0*\edgeDist,0*\edgeDist) circle (2pt) node at (0*\edgeDist,0*\edgeDist-\putLeft) {$2$};
\filldraw[black] (0*\edgeDist,1*\edgeDist) circle (2pt) node at (0*\edgeDist,1*\edgeDist+\putLeft) {$1$};
\filldraw[black] (1*\edgeDist,0*\edgeDist) circle (2pt) node at (1*\edgeDist,0*\edgeDist-\putLeft) {$4$};
\filldraw[black] (1*\edgeDist,1*\edgeDist) circle (2pt) node at (1*\edgeDist,1*\edgeDist+\putLeft) {$3$};
\filldraw[black] (2*\edgeDist,0*\edgeDist) circle (2pt) node at (2*\edgeDist,0*\edgeDist-\putLeft) {$6$};
\filldraw[black] (2*\edgeDist,1*\edgeDist) circle (2pt) node at (2*\edgeDist,1*\edgeDist+\putLeft) {$5$};
%
\draw (0*\edgeDist,0*\edgeDist) to (0*\edgeDist,1*\edgeDist); 
\draw (0*\edgeDist,0*\edgeDist) to (1*\edgeDist,1*\edgeDist); 
\draw (1*\edgeDist,0*\edgeDist) to (1*\edgeDist,1*\edgeDist); 
\draw (2*\edgeDist,1*\edgeDist) to (0*\edgeDist,0*\edgeDist); 
\draw (2*\edgeDist,1*\edgeDist) to (1*\edgeDist,0*\edgeDist); 
\draw (2*\edgeDist,1*\edgeDist) to (2*\edgeDist,0*\edgeDist); 
\end{tikzpicture}
\caption{$F_{3}$}
\end{subfigure}
\hfill
\begin{subfigure}{.4\textwidth}
\centering
\begin{tikzpicture}
\filldraw[black] (0*\edgeDist,0*\edgeDist) circle (2pt) node at (0*\edgeDist,0*\edgeDist-\putLeft) {$2$};
\filldraw[black] (0*\edgeDist,1*\edgeDist) circle (2pt) node at (0*\edgeDist,1*\edgeDist+\putLeft) {$1$};
\filldraw[black] (1*\edgeDist,0*\edgeDist) circle (2pt) node at (1*\edgeDist,0*\edgeDist-\putLeft) {$4$};
\filldraw[black] (1*\edgeDist,1*\edgeDist) circle (2pt) node at (1*\edgeDist,1*\edgeDist+\putLeft) {$3$};
\filldraw[black] (2*\edgeDist,0*\edgeDist) circle (2pt) node at (2*\edgeDist,0*\edgeDist-\putLeft) {$6$};
\filldraw[black] (2*\edgeDist,1*\edgeDist) circle (2pt) node at (2*\edgeDist,1*\edgeDist+\putLeft) {$5$};
\filldraw[black] (3*\edgeDist,0*\edgeDist) circle (2pt) node at (3*\edgeDist,0*\edgeDist-\putLeft) {$8$};
\filldraw[black] (3*\edgeDist,1*\edgeDist) circle (2pt) node at (3*\edgeDist,1*\edgeDist+\putLeft) {$7$};
%
\draw (0*\edgeDist,0*\edgeDist) to (0*\edgeDist,1*\edgeDist); 
\draw (0*\edgeDist,0*\edgeDist) to (1*\edgeDist,1*\edgeDist); 
\draw (1*\edgeDist,0*\edgeDist) to (1*\edgeDist,1*\edgeDist); 
\draw (2*\edgeDist,1*\edgeDist) to (0*\edgeDist,0*\edgeDist); 
\draw (2*\edgeDist,1*\edgeDist) to (1*\edgeDist,0*\edgeDist); 
\draw (2*\edgeDist,1*\edgeDist) to (2*\edgeDist,0*\edgeDist); 
\draw (3*\edgeDist,1*\edgeDist) to (0*\edgeDist,0*\edgeDist); 
\draw (3*\edgeDist,1*\edgeDist) to (1*\edgeDist,0*\edgeDist); 
\draw (3*\edgeDist,1*\edgeDist) to (2*\edgeDist,0*\edgeDist); 
\draw (3*\edgeDist,1*\edgeDist) to (3*\edgeDist,0*\edgeDist); 
\end{tikzpicture}
\caption{$F_{4}$}
\end{subfigure}
\caption{$F_{m}$, $m \leq 4$}
\label{fig_Fm}
\end{figure}

\bthm[{\cite[Theorem 6.1]{bolognini2018binomial}}]
\label{thm_char_bipartite_CM_graphs}
Let $G$ be a connected bipartite graph. The following properties are equivalent:
\begin{enumerate}
    \item $J_{G}$ is Cohen--Macaulay,
    \item There exists $s \geq 1$ such that $G = G_{1} \ast \cdots \ast G_{s}$, where $G_{i} = F_{n_{i}}$ or ${G_{i} = F_{m_{i,1}} \circ \cdots \circ F_{m_{i,t_{i}}}}$, for some $n_{i} \geq 1$ and $m_{i,j} \geq 3$ for each $j = 1,\ldots,t_{i}$.
\end{enumerate}
\ethm

With the decomposition of $G$ as in Theorem \ref{thm_char_bipartite_CM_graphs}, define the following: $A = \{i \in [s] \mid G_{i} = F_{n_{i}}, n_{i} \geq 2\}$, $B = \{i \in [s] \mid G_{i} = F_{n_{i}}, n_{i} = 1\}$, and $C = \{i \in [s] \mid G_{i} = F_{m_{i,1}} \circ \cdots F_{m_{i,t_{i}}}, t_{i} \geq 2\}$. For each $i \in C$, let $C_{i} = \{j \in \{2,\ldots,t_{i}-1\} \mid m_{i,j} \geq 4\} \cup \{1,t_{i}\}$ and $C_{i}^{'} = \{j \in \{2,\ldots,t_{i}-1\} \mid m_{i,j} = 3\}$. Set $\alpha = \card{A} + \sum_{i \in C} \card{C_{i}}$ and $\beta = \card{B} + \sum_{i\in C} \card{C_{i}^{'}}$.

\bthm[{\cite[Theorem 4.7]{jayanthan2019regularityCM&BipartiteGraphs}}]
\label{thm_reg_CM_bipartite}
Let $G = G_{1}\ast \cdots \ast G_{s}$ be a Cohen--Macaulay connected bipartite graph. Let $\alpha$ and $\beta$ be defined as above, then $\reg(S/J_{G}) = 3\alpha + \beta$.
\ethm

\bcor
\label{cor_nu_equals_reg_CM_bipartite_graphs}
Let $G = G_{1}\ast \cdots \ast G_{s}$ be a Cohen--Macaulay connected bipartite graph. Then, $\nu(G) = \reg(S/J_{G})$.
\ecor

\bproof
For convenience of the proof we assume that the vertices of each $G_{i}$ are labeled by the integers $1$ through $\card{V(G_{i})}$.
For $1 \leq i \leq s$, we construct a path $P_{i}$ inside $G_{i}$ as follows:
\begin{enumerate}
    \item If $i \in A$, let $P_{i}$ be the path $[1,2,3,4]$ on $G_{i}$,
    \item If $i \in B$, let $P_{i}$ be the path $[1,2]$ on $G_{i}$,
    \item Suppose $i \in C$ and that $G_{i} = F_{m_{i,1}} \circ \cdots F_{m_{i,t_{i}}}, t_{i} \geq 2$. For $1 \leq j \leq m_{i,t_{i}}$ construct a path $P_{i,j}$ inside $F_{m_{i,j}}$ as follows:
    \begin{enumerate}
        \item If $j = 1$, let $P_{i,j}$ be the path $[1,2,3,4]$ on $F_{m_{i,1}}$.
        \item If $j = t_{i}$, let $P_{i,j}$ be the path $[2\cdot m_{i,t_{i}}-3,2\cdot m_{i,t_{i}}-2,2\cdot m_{i,t_{i}}-1,2\cdot m_{i,t_{i}}]$ on $F_{m_{i,t_{i}}}$.
        \item If $j \in C_{i} \setminus \{1,t_{i}\}$, let $P_{i,j}$ be the path $[3,4,5,6]$ on $F_{m_{i,j}}$.
        \item If $j \in C_{i}^{'}$, let $P_{i,j}$ be the path $[3,4]$ on $F_{m_{i,j}}$.
    \end{enumerate}
    Put $P_{i} = \bigcup_{j=1}^{m_{i,t_{i}}} P_{i,j}$.
\end{enumerate}
We define $\ul{P}$ as $\bigcup_{i=1}^{s} P_{i}$ and $P_{\ind}$ as the induced subgraph of $G$ on $V(\ul{P})$. The construction of the $P_{i}$ yields that for $1 \leq i \leq s-1$, either $P_{i}$ and $P_{i+1}$ share a terminal vertex or there is no induced edge between $P_{i}$ and $P_{i+1}$. The construction of $G$ via Definition \ref{defn_operations_bipartite_graphs} implies that there are no edges between $V(G_{i})$ and $V(G_{j})$ whenever $\card{i-j} > 1$. Hence, in particular, $P_{i}$ and $P_{j}$ have no induced edge whenever $\card{i - j} > 1$. It follows that $P_{\ind}$ is a disjoint union of induced paths. Hence, in particular, $\ul{P}$ is DOIP. Thus, 
\begin{align*}
    3\alpha + \beta 
    &= \sum_{i=1}^{s} \card{E(P_{i})} & & (\text{by construction of the } P_{i}) \\
    &\leq \nu(G) & & (\ul{P} \text{ is DOIP}) \\
    &\leq \reg(S/J_{G}) & & (\text{Theorem } \ref{thm_nu_leq_reg}) \\
    &= 3\alpha + \beta. & & (\text{Theorem } \ref{thm_reg_CM_bipartite})
\end{align*}
\eproof

\bex
\label{ex_nu_equals_reg_fomula}
We consider the Cohen--Macaulay bipartite graph $G = G_{1} \ast G_{2} \ast G_{3}$ where $G_{1} = (F_{3}\circ F_{4} \circ F_{3} \circ F_{4})$, $G_{2} = F_{1}$, and $G_{3} = F_{4}$. We illustrate the construction of $\ul{P}$ in the proof of Corollary \ref{cor_nu_equals_reg_CM_bipartite_graphs} via Figure \ref{fig_total_construction}. The graph $G_{1}$ is constructed by identifying:
\begin{enumerate}
    \item vertex $5$ of subfigure \ref{fig_g1_1} with vertex $2$ of subfigure \ref{fig_g1_2},
    \item vertex $7$ of subfigure \ref{fig_g1_2} with vertex $2$ of subfigure \ref{fig_g1_3},
    \item vertex $5$ of subfigure \ref{fig_g1_3} with vertex $2$ of subfigure \ref{fig_g1_4}.
\end{enumerate}
The graph $G$ is constructed from $G_{1}$, $G_{2}$, and $G_{3}$ by identifying:
\begin{enumerate}
    \item vertex $8$ of subfigure \ref{fig_g1_4} with vertex $1$ of Figure \ref{fig_g2},
    \item vertex $2$ of Figure \ref{fig_g2} with vertex $1$ of Figure \ref{fig_g3}.
\end{enumerate}
The proof of Corollary \ref{cor_nu_equals_reg_CM_bipartite_graphs} yields that:
\begin{enumerate}
    \item for $G_{1}$, we have that $P_{1,1} = [1,2,3,4]$, $P_{1,2} = [3,4,5,6]$, $P_{1,3} = [3,4]$, and $P_{1,4} = [5,6,7,8]$. Then, $P_{1}$ is the disjoint union of these paths in $G_{1}$.
    \item for $G_{2}$, we have that $P_{2} = [1,2]$, and
    \item for $G_{3}$, we have that $P_{3} = [1,2,3,4]$.
\end{enumerate}
It follows that $\ul{P} := \bigcup_{i=1}^{\ell} P_{i}$ consists of four vertex-disjoint induced paths. (In the gluing step, terminal vertices of $P_{1,4}$, $P_{2}$, and $P_{3}$ coincide. These paths concatenate in the construction of $G$.) As mentioned in the proof, all four of these vertex-disjoint induced paths contain no induced edges between themselves.
\begin{figure}
\begin{subfigure}{.14\textwidth}
\centering
\begin{tikzpicture}
\filldraw[black] (0*\edgeDist,0*\edgeDist) circle (2pt) node at (0*\edgeDist,0*\edgeDist-\putLeft) {$2$};
\filldraw[black] (0*\edgeDist,1*\edgeDist) circle (2pt) node at (0*\edgeDist,1*\edgeDist+\putLeft) {$1$};
\filldraw[black] (1*\edgeDist,0*\edgeDist) circle (2pt) node at (1*\edgeDist,0*\edgeDist-\putLeft) {$4$};
\filldraw[black] (1*\edgeDist,1*\edgeDist) circle (2pt) node at (1*\edgeDist,1*\edgeDist+\putLeft) {$3$};
\filldraw[black] (2*\edgeDist,1*\edgeDist) circle (2pt) node at (2*\edgeDist,1*\edgeDist+\putLeft) {$5$};
%
\draw[red] (0*\edgeDist,0*\edgeDist) to (0*\edgeDist,1*\edgeDist); 
\draw[red] (0*\edgeDist,0*\edgeDist) to (1*\edgeDist,1*\edgeDist); 
\draw[red] (1*\edgeDist,0*\edgeDist) to (1*\edgeDist,1*\edgeDist); 
\draw (2*\edgeDist,1*\edgeDist) to (0*\edgeDist,0*\edgeDist); 
\draw (2*\edgeDist,1*\edgeDist) to (1*\edgeDist,0*\edgeDist); 
\end{tikzpicture}
\caption{$F_{3} \setminus \{6\}$}
\label{fig_g1_1}
\end{subfigure}
\hfill
\begin{subfigure}{.245\textwidth}
\centering
\begin{tikzpicture}
\filldraw[black] (0*\edgeDist,0*\edgeDist) circle (2pt) node at (0*\edgeDist,0*\edgeDist-\putLeft) {$2$};
\filldraw[black] (1*\edgeDist,0*\edgeDist) circle (2pt) node at (1*\edgeDist,0*\edgeDist-\putLeft) {$4$};
\filldraw[black] (1*\edgeDist,1*\edgeDist) circle (2pt) node at (1*\edgeDist,1*\edgeDist+\putLeft) {$3$};
\filldraw[black] (2*\edgeDist,0*\edgeDist) circle (2pt) node at (2*\edgeDist,0*\edgeDist-\putLeft) {$6$};
\filldraw[black] (2*\edgeDist,1*\edgeDist) circle (2pt) node at (2*\edgeDist,1*\edgeDist+\putLeft) {$5$};
\filldraw[black] (3*\edgeDist,1*\edgeDist) circle (2pt) node at (3*\edgeDist,1*\edgeDist+\putLeft) {$7$};
%
\draw (0*\edgeDist,0*\edgeDist) to (1*\edgeDist,1*\edgeDist); 
\draw[red] (1*\edgeDist,0*\edgeDist) to (1*\edgeDist,1*\edgeDist); 
\draw (2*\edgeDist,1*\edgeDist) to (0*\edgeDist,0*\edgeDist); 
\draw[red] (2*\edgeDist,1*\edgeDist) to (1*\edgeDist,0*\edgeDist); 
\draw[red] (2*\edgeDist,1*\edgeDist) to (2*\edgeDist,0*\edgeDist); 
\draw (3*\edgeDist,1*\edgeDist) to (0*\edgeDist,0*\edgeDist); 
\draw (3*\edgeDist,1*\edgeDist) to (1*\edgeDist,0*\edgeDist); 
\draw (3*\edgeDist,1*\edgeDist) to (2*\edgeDist,0*\edgeDist); 
\end{tikzpicture}
\caption{$F_{4} \setminus \{1,8\}$}
\label{fig_g1_2}
\end{subfigure}
\hfill 
\begin{subfigure}{.15\textwidth}
\centering
\begin{tikzpicture}
\filldraw[black] (0*\edgeDist,0*\edgeDist) circle (2pt) node at (0*\edgeDist,0*\edgeDist-\putLeft) {$2$};
\filldraw[black] (1*\edgeDist,0*\edgeDist) circle (2pt) node at (1*\edgeDist,0*\edgeDist-\putLeft) {$4$};
\filldraw[black] (1*\edgeDist,1*\edgeDist) circle (2pt) node at (1*\edgeDist,1*\edgeDist+\putLeft) {$3$};
\filldraw[black] (2*\edgeDist,1*\edgeDist) circle (2pt) node at (2*\edgeDist,1*\edgeDist+\putLeft) {$5$};
%
\draw (0*\edgeDist,0*\edgeDist) to (1*\edgeDist,1*\edgeDist); 
\draw[red] (1*\edgeDist,0*\edgeDist) to (1*\edgeDist,1*\edgeDist); 
\draw (2*\edgeDist,1*\edgeDist) to (0*\edgeDist,0*\edgeDist); 
\draw (2*\edgeDist,1*\edgeDist) to (1*\edgeDist,0*\edgeDist); 
\end{tikzpicture}
\caption{$F_{3} \setminus \{1,6\}$}
\label{fig_g1_3}
\end{subfigure}
\hfill
\begin{subfigure}{.225\textwidth}
\centering
\begin{tikzpicture}
\filldraw[black] (0*\edgeDist,0*\edgeDist) circle (2pt) node at (0*\edgeDist,0*\edgeDist-\putLeft) {$2$};
\filldraw[black] (1*\edgeDist,0*\edgeDist) circle (2pt) node at (1*\edgeDist,0*\edgeDist-\putLeft) {$4$};
\filldraw[black] (1*\edgeDist,1*\edgeDist) circle (2pt) node at (1*\edgeDist,1*\edgeDist+\putLeft) {$3$};
\filldraw[black] (2*\edgeDist,0*\edgeDist) circle (2pt) node at (2*\edgeDist,0*\edgeDist-\putLeft) {$6$};
\filldraw[black] (2*\edgeDist,1*\edgeDist) circle (2pt) node at (2*\edgeDist,1*\edgeDist+\putLeft) {$5$};
\filldraw[black] (3*\edgeDist,1*\edgeDist) circle (2pt) node at (3*\edgeDist,1*\edgeDist+\putLeft) {$7$};
\filldraw[black] (3*\edgeDist,0*\edgeDist) circle (2pt) node at (3*\edgeDist,0*\edgeDist-\putLeft) {$8$};
%
\draw (0*\edgeDist,0*\edgeDist) to (1*\edgeDist,1*\edgeDist); 
\draw (1*\edgeDist,0*\edgeDist) to (1*\edgeDist,1*\edgeDist); 
\draw (2*\edgeDist,1*\edgeDist) to (0*\edgeDist,0*\edgeDist); 
\draw (2*\edgeDist,1*\edgeDist) to (1*\edgeDist,0*\edgeDist); 
\draw[red] (2*\edgeDist,1*\edgeDist) to (2*\edgeDist,0*\edgeDist); 
\draw (3*\edgeDist,1*\edgeDist) to (0*\edgeDist,0*\edgeDist); 
\draw (3*\edgeDist,1*\edgeDist) to (1*\edgeDist,0*\edgeDist); 
\draw[red] (3*\edgeDist,1*\edgeDist) to (2*\edgeDist,0*\edgeDist); 
\draw[red] (3*\edgeDist,1*\edgeDist) to (3*\edgeDist,0*\edgeDist); 
\end{tikzpicture}
\caption{$F_{4} \setminus \{1\}$}
\label{fig_g1_4}
\end{subfigure}
\caption*{The components of $G_{1}$}

\hfill
\begin{subfigure}{.4\textwidth} 
\centering
\begin{tikzpicture}
\filldraw[black] (0*\edgeDist,0*\edgeDist) circle (2pt) node at (0*\edgeDist,0*\edgeDist-\putLeft) {$2$};
\filldraw[black] (0*\edgeDist,1*\edgeDist) circle (2pt) node at (0*\edgeDist,1*\edgeDist+\putLeft) {$1$};
%
\draw[red] (0*\edgeDist,0*\edgeDist) to (0*\edgeDist,1*\edgeDist); 
\end{tikzpicture}
\caption{$G_{2} = F_{1}$}
\label{fig_g2}
\end{subfigure}
\hfill
%
\begin{subfigure}{.4\textwidth}
\centering
\begin{tikzpicture}
\filldraw[black] (0*\edgeDist,0*\edgeDist) circle (2pt) node at (0*\edgeDist,0*\edgeDist-\putLeft) {$2$};
\filldraw[black] (0*\edgeDist,1*\edgeDist) circle (2pt) node at (0*\edgeDist,1*\edgeDist+\putLeft) {$1$};
\filldraw[black] (1*\edgeDist,0*\edgeDist) circle (2pt) node at (1*\edgeDist,0*\edgeDist-\putLeft) {$4$};
\filldraw[black] (1*\edgeDist,1*\edgeDist) circle (2pt) node at (1*\edgeDist,1*\edgeDist+\putLeft) {$3$};
\filldraw[black] (2*\edgeDist,0*\edgeDist) circle (2pt) node at (2*\edgeDist,0*\edgeDist-\putLeft) {$6$};
\filldraw[black] (2*\edgeDist,1*\edgeDist) circle (2pt) node at (2*\edgeDist,1*\edgeDist+\putLeft) {$5$};
\filldraw[black] (3*\edgeDist,0*\edgeDist) circle (2pt) node at (3*\edgeDist,0*\edgeDist-\putLeft) {$8$};
\filldraw[black] (3*\edgeDist,1*\edgeDist) circle (2pt) node at (3*\edgeDist,1*\edgeDist+\putLeft) {$7$};
%
\draw[red] (0*\edgeDist,0*\edgeDist) to (0*\edgeDist,1*\edgeDist); 
\draw[red] (0*\edgeDist,0*\edgeDist) to (1*\edgeDist,1*\edgeDist); 
\draw[red] (1*\edgeDist,0*\edgeDist) to (1*\edgeDist,1*\edgeDist); 
\draw (2*\edgeDist,1*\edgeDist) to (0*\edgeDist,0*\edgeDist); 
\draw (2*\edgeDist,1*\edgeDist) to (1*\edgeDist,0*\edgeDist); 
\draw (2*\edgeDist,1*\edgeDist) to (2*\edgeDist,0*\edgeDist); 
\draw (3*\edgeDist,1*\edgeDist) to (0*\edgeDist,0*\edgeDist); 
\draw (3*\edgeDist,1*\edgeDist) to (1*\edgeDist,0*\edgeDist); 
\draw (3*\edgeDist,1*\edgeDist) to (2*\edgeDist,0*\edgeDist); 
\draw (3*\edgeDist,1*\edgeDist) to (3*\edgeDist,0*\edgeDist); 
\end{tikzpicture}
\caption{$G_{3} = F_{4}$}
\label{fig_g3}
\end{subfigure}
\caption{Construction of $G$ from $G_{1}$, $G_{2}$, and $G_{3}$}
\label{fig_total_construction}
\end{figure}
\eex 

\section{Preliminaries on Block Graphs}
\label{sec_prelims_block_graphs}

In this section, we recall the definition and some elementary properties of block graphs. All of the results in this section are well-known to experts. However, lacking a reference for these results, we present their proofs for completeness.

Recall the definition of a block graph.
\bdefn[Block Graph]
\label{defn_block_graph}
A  graph $G$ is \textbf{biconnected} if $G$ is connected and $G \setminus v$ is connected for every $v \in V(G)$. $G$ is a \textbf{clique} or a \textbf{complete graph} if for every pair of distinct vertices $v$ and $w$ in $V(G)$, $\{v,w\} \in E(G)$. A subgraph $B$ of $G$ is a \textbf{block} of $G$ if $B$ is a maximal biconnected component of $G$ with respect to inclusion. A graph $G$ is a \textbf{block graph} if every block of $G$ is a complete graph. Block graphs are also referred to in the literature as a \textbf{tree of cliques}.
\edefn 

We recall the following well-known properties of block graphs.

\blem
\label{lem_induced_subgraph_block_graph_block_graph}
Let $G$ be a block graph. If $H$ is an induced subgraph of $G$, then $H$ is a block graph.
\elem

\blem
\label{lem_shortest_path_block_graph_induced}
Let $G$ be a connected block graph. Let $v$ and $w$ be distinct vertices of $G$. There is a unique shortest path in $G$ connecting $v$ and $w$, and this shortest path is an induced path in $G$.
\elem 

\blem
\label{lem_block_graphs_do_not_have_certain_cycles}
Let $G$ be a block graph. Let $Q$ be an induced path in $G$. Then, every cycle of $G$ contains no more than one edge of $Q$.
\elem

\bproof
Suppose by contradiction that there exists a cycle, $C$ of $G$, which contains two or more edges of $Q$. Since $G$ is a block graph, $C$ induces a complete subgraph of $G$ which would contradict $Q$ being induced.
\eproof

\blem
\label{lem_decomp_intersecting_induced_paths}
Let $Q = [v_{1},\ldots,v_{r}]$ and $R = [w_{1},\ldots,w_{s}]$ be induced paths of a block graph which intersect non-trivially. Then, $Q \cap R$ is connected. Moreover, there is a unique decomposition $Q = Q_{1} \ast \gamma \ast Q_{2}$ and $R = R_{1} \ast \gamma \ast R_{2}$ where $Q_{i}$ (respectively $R_{i}$) is possibly a singleton vertex of $Q$ (respectively $R$) for $i = 1,2$, and $\gamma = Q \cap R$.
\elem

\bproof
Let $v$ and $w$ be vertices of $Q \cap R$. Let $[v,w]\mid_{Q}$ and $[v,w]\mid_{R}$ be the subpaths of $Q$ and $R$ from $v$ to $w$. These are induced subpaths of $Q$ and $R$, respectively, and hence are induced paths of $G$, by Lemma \ref{lem_transitivity_of_being_induced}. Lemma \ref{lem_shortest_path_block_graph_induced} implies that these are the same path. Thus, $Q \cap R$ is connected.
\eproof

\blem
\label{lem_decomp_induced_path_intersecting_two_others}
Let $R_{1} = [w_{1},\ldots,w_{s}]$ and $R_{2} = [w_{1}^{'},\ldots,w_{t}^{'}]$ be disjoint induced paths of a block graph, and $Q = [v_{1},\ldots,v_{r}]$ be an induced path of a block graph. Suppose that $v_{1} = w_{1}$ and $v_{r} = w_{t}^{'}$. Then, there is a unique decomposition of $Q$ as $Q_{1} \ast \gamma \ast Q_{2}$ where $Q_{i}$ is a subpath of $R_{i}$ for $i = 1,2$, and $\gamma \cap R_{i}$ is a vertex for $i = 1,2$.
\elem

\bproof
Apply Lemma \ref{lem_decomp_intersecting_induced_paths} to $R_{1}$ and $Q$. Because $R_{1} \cap Q$ contains $v_{1}$, $Q = \gamma^{'} \ast Q_{2}^{'}$, where $\gamma^{'} := R_{1} \cap Q$ and $Q_{2}^{'}$ intersects $R_{1}$ at a vertex. Because $R_{1}$ and $R_{2}$ are disjoint, $Q_{2}^{'}$ contains $w_{t}^{'}$. Apply Lemma \ref{lem_decomp_intersecting_induced_paths} to $R_{2}$ and $Q_{2}^{'}$. Then, $Q_{2}^{'} = Q_{1}^{''} \ast \gamma^{''}$ where $\gamma^{''} := R_{2} \cap Q_{2}^{'}$ and $Q_{1}^{''}$ intersects $R_{2}$ at a vertex. Put $Q_{1} := \gamma^{'}$, $\gamma := Q_{1}^{''}$, and $Q_{2} := \gamma^{''}$. Then, $Q = Q_{1} \ast \gamma \ast Q_{2}$.
\eproof

\section{Combinatorial Characterization of DOIP Paths in Block Graphs} \label{sec_block_graph}

In this section, we give several equivalent combinatorial formulations for vertex-disjoint paths of a block graph $G$ to be DOIP, which are in terms of forbidden subgraphs of $P_{\ind}$. 

\subsection{Forbidden Subgraphs}
Let $G$ be a block graph, and let $\ul{P} := P_{1},\ldots,P_{\ell}$ be vertex-disjoint induced oriented paths of $G$.

\bdefn 
\label{defn_strand}
Let $Q$ be an oriented induced path of $P_{\ind}$ with orientation $\phi_{Q}$. We denote by $Q_{0} := \phi_{Q}(1)$ and $Q_{r} := \phi_{Q}(2)$. We say that $Q$ is a \textbf{strand} of $P_{\ind}$ from $P_{i}$ to $P_{j}$, with $i \neq j$, if:
\begin{enumerate}
    \item $Q_{0} \neq \phi_{i}(2)$,
    \item $Q_{r} \neq \phi_{j}(1)$,
    \item $V(Q) \cap V(P_{i}) = \{Q_{0}\}$,
    \item $V(Q) \cap V(P_{j}) = \{Q_{r}\}$,
    \item $Q$ does not contain $P_{k}$ for any $1 \leq k \leq \ell$,
\end{enumerate}
We say that $Q$ is an \textbf{internal strand} if $Q$ is a strand from $P_{i}$ to $P_{j}$ and, in addition, $Q_{0} \in P_{i}^{\circ}$ and $Q_{r} \in P_{j}^{\circ}$. 
\edefn 

\brem 
We observe that conditions (1) and (2) of Definition \ref{defn_strand} will be automatically satisfied whenever $H$ is an internal strand.
\erem 

\brem 
The motivation for Definition \ref{defn_strand} is that it will allow us to relate strands and arcs of $K_{P_{\ind}}$.
\erem

\bdefn 
\label{defn_fork}
Let $i$ and $j$ be distinct integers belonging to $[\ell]$, $a$ a vertex belonging to $V(P_{i}) \setminus \{ \phi_{i}(2)\}$, and $b$ and $c$ distinct vertices belonging to $V(P_{j})$. Let $Q$ be an oriented induced path with orientation $\phi_{Q}$ satisfying:
\begin{enumerate}
    \item $\phi_{Q}(1) = a$,
    \item $V(Q) \cap V(P_{i}) = \{a\}$,
    \item $V(Q) \cap V(P_{j}) = \varnothing$,
    \item $\phi_{Q}(2)$ is adjacent to $b$ and $c$ in $P_{\ind}$,
    \item $Q$ does not contain $P_{k}$ for any $1 \leq k \leq \ell$.
\end{enumerate}
We define the subgraph $H$ of $P_{\ind}$ as follows: 
\begin{align*}
    V(H) &:= V(Q) \cup \{b,c\} \\
    E(H) &:= E(Q) \cup \{ \{\phi_{Q}(2),b\}, \{\phi_{Q}(2),c\} \}.
\end{align*}
We say that $H$ is a \textbf{fork} of $P_{\ind}$ from $P_{i}$ to $P_{j}$. We say that $H$ is an \textbf{internal fork} of $P_{\ind}$ from $P_{i}$ to $P_{j}$ if, in addition, $a \in P_{i}^{\circ}$. 

For a vertex $v \in H$, we say that $v$ is a \textbf{terminal vertex} of $H$ if $v \in \{a,b,c\}$; otherwise, we say that $v$ is an \textbf{internal vertex}.
\edefn 

\brem 
Lemma \ref{lem_block_graphs_do_not_have_certain_cycles} implies that $\{b,c\} \in E(P_{\ind})$.
\erem 

\brem 
We observe that the requirement $a \in V(P_{i}) \setminus \{ \phi_{i}(2)\}$ will be automatically satisfied whenever $H$ is an internal fork.
\erem 

\bdefn 
\label{defn_double_fork}
Let $i$ and $j$ be distinct integers belonging to $[\ell]$, $a$ and $b$ distinct vertices of $V(P_{i})$, and $c$ and $d$ distinct vertices belonging to $V(P_{j})$. Let $Q$ be an oriented induced path with orientation $\phi_{Q}$ satisfying:
\begin{enumerate}
    \item $\phi_{Q}(1)$ is adjacent to $a$ and $b$ in $P_{\ind}$,
    \item $\phi_{Q}(2)$ is adjacent to $c$ and $d$ in $P_{\ind}$,
    \item $V(Q) \cap V(P_{i}) = \varnothing$,
    \item $V(Q) \cap V(P_{j}) = \varnothing$,
    \item $Q$ does not contain $P_{k}$ for any $1 \leq k \leq \ell$.
\end{enumerate}
We define the subgraph $H$ of $P_{\ind}$ as follows: 
\begin{align*}
    V(H) &:= V(Q) \cup \{a,b,c,d\} \\
    E(H) &:= E(Q) \cup \{ \{\phi_{Q}(1),a\}, \{\phi_{Q}(1),b\}, \{\phi_{Q}(2),c\}, \{\phi_{Q}(2),d\} \}.
\end{align*}
We say that $H$ is a \textbf{double fork} of $P_{\ind}$ from $P_{i}$ to $P_{j}$.

For a vertex $v \in H$, we say that $v$ is a \textbf{terminal vertex} of $H$ if $v \in \{a,b,c,d\}$; otherwise, we say that $v$ is an \textbf{internal vertex}.
\edefn 

\bdefn 
\label{defn_complete_ladder}
Let $i$ and $j$ be distinct integers belonging to $[\ell]$, $a$ and $b$ distinct vertices of $V(P_{i})$, and $c$ and $d$ distinct vertices belonging to $V(P_{j})$. Let $H$ denote the induced subgraph on $\{a,b,c,d\}$. We say that $H$ is a \textbf{complete ladder} if $H$ is a complete graph.
\edefn

\brem 
We observe that in Definition \ref{defn_complete_ladder}, a complete ladder is $K_{4}$, the complete graph on $4$ vertices. However, not every $K_{4}$ in $P_{\ind}$ is a complete ladder. For example, a $K_{4}$ whose vertices are terminal vertices of distinct $P_{i}$ would not be a complete ladder, nor would it realize an internal strand, an internal fork, or a double fork.
\erem 

\bex 
\label{ex_forbidden_subgraphs}
Let $G$ be the graph in Figure \ref{fig_forbidden_paths}. We observe that $G$ is a block graph. We consider the vertex-disjoint induced paths
\begin{align*}
    P_{1} &:= [1,2,3] & P_{2} &:= [4,5,6] & P_{3} &:= [7,8] & & \\
    P_{4} &:= [9,10,11] & P_{5} &:= [12,13] & P_{6} &:= [14,15,16] & P_{7} &:= [17,18] \\
    P_{8} &:= [19,20] & P_{9} &:= [21,22,23] & P_{10} &:= [24,25,26] & P_{11} &:= [27,28,29].
\end{align*}
For $1 \leq i \leq 11$, we define the orientation $\phi_{i}$ such that $\phi_{i}(1) < \phi_{i}(2)$. 
We present some subgraphs of $G$ that illustrate Definitions \ref{defn_strand}, \ref{defn_fork}, \ref{defn_double_fork},  and \ref{defn_complete_ladder}.
\begin{enumerate}
    \item (Internal) Strand:
    \begin{enumerate}
        \item $Q = [2,5]$ is an (internal) strand from $P_{1}$ to $P_{2}$, $\phi_{Q}(1) = 2$, and $\phi_{Q}(2) = 5$,
        \item $[5,10,11,13,15]$ is an (internal) strand from $P_{2}$ to $P_{6}$. 
        \item $[5,2]$ is a strand from $P_{2}$ to $P_{1}$,
        \item $[17,20]$ is a strand from $P_{7}$ to $P_{8}$,
        \item $[2,5,10,11,13,15,18,19,21,25,28]$ is an internal strand from $P_{1}$ to $P_{11}$.
    \end{enumerate}
    \item (Internal) Fork: (we just list the vertices of the fork)
    \begin{enumerate}
        \item $\{7,8,9\}$ is a fork from $P_{4}$ to $P_{3}$, $Q = [9]$ is the singleton path,
        \item $\{7,8,9,10,5\}$ is an internal fork from $P_{2}$ to $P_{3}$, $Q = [5,10,9]$, $\phi_{Q}(1) = 5$, $\phi_{Q}(2) = 9$,
        \item $\{17,18,19,21\}$ is a fork from $P_{9}$ to $P_{7}$,
        \item $\{22,25,27,28\}$ is an internal fork from $P_{9}$ to $P_{11}$.
    \end{enumerate}
    \item Double Fork: (we just list the vertices of the double fork)
    \begin{enumerate}
        \item $\{21,22,25,27,28\}$ is a double fork, $Q = [25]$ is a singleton path,
        \item $\{17,18,19,21,25,27,28\}$ is a double fork, $Q = [19,21,25]$, $\phi_{Q}(1) = 19$, $\phi_{Q}(2) = 25$.
    \end{enumerate}
    \item Complete Ladder: (we just list the vertices of the double fork)
    \begin{enumerate}
        \item $\{17,18,19,20\}$ is a complete ladder.
    \end{enumerate}
\end{enumerate}
We next present examples of subgraphs of $G$ that do not satisfy the requirements of Definitions \ref{defn_strand}, \ref{defn_fork}, \ref{defn_double_fork}, and \ref{defn_complete_ladder}.
\begin{enumerate}
    \item Not strands:
    \begin{enumerate}
        \item $[8,9,10,5]$ is not a strand from $P_{3}$ to $P_{2}$ because $\phi_{Q}(1) = 8 = \phi_{3}(2)$,
        \item $[7,9]$ is not a strand from $P_{3}$ to $P_{4}$ because $\phi_{Q}(2) = 9 = \phi_{4}(1)$,
        \item $[8,9]$ is not a strand from $P_{3}$ to $P_{4}$ (same reason), but it is a strand from $P_{4}$ to $P_{3}$,
        \item $[1,2,5]$ is not a strand from $P_{1}$ to $P_{2}$ because $Q$ violates condition (3) of Definition \ref{defn_strand},
        \item $[7,9,10,11,13]$ is not a strand from $P_{3}$ to $P_{5}$ because $Q$ contains $P_{4}$.
    \end{enumerate}
    \item Not forks: (we just list the vertices under consideration)
    \begin{enumerate}
        \item $\{12,13,15,18\}$ is not a fork because $\phi_{Q}(1) = 18 = \phi_{7}(2)$,
        \item $\{2,3,5,4\}$ is not a fork because it violates condition (2) of Definition \ref{defn_fork},
        \item $\{7,8,9,10,11,13,15\}$ is not a fork from $P_{6}$ to $P_{3}$ because $Q$ contains $P_{4}$.
    \end{enumerate}
    \item Not double forks: (we just list the vertices under consideration)
    \begin{enumerate}
        \item $\{7,8,9,10,11,13,15,18,19,20\}$ is not a double fork because $Q = [9,10,11,13,14,18]$ contains $P_{4}$.
    \end{enumerate}
\end{enumerate}
\eex 

\def\putLeft{.5*.9}
\def\edgeDist{1.75*.9}
\def\iii{3}
\def\iv{4}
\def\v{5}
\def\vi{6}
\def\vii{7}
\def\viii{8}
\def\mi{-1}
\begin{figure}
\begin{center}
\begin{tikzpicture}
\filldraw[black] (-2*\edgeDist,1*\edgeDist) circle (2pt) node at (-2*\edgeDist,1*\edgeDist+\putLeft) {$1$};
\filldraw[black] (-2*\edgeDist,0*\edgeDist) circle (2pt) node at (-2*\edgeDist+\putLeft,0*\edgeDist+\putLeft) {$2$};
\filldraw[black] (-2*\edgeDist,-1*\edgeDist) circle (2pt) node at (-2*\edgeDist+\putLeft,-1*\edgeDist) {$3$};

\filldraw[black] (-1*\edgeDist,1*\edgeDist) circle (2pt) node at (-1*\edgeDist,1*\edgeDist+\putLeft) {$4$};
\filldraw[black] (-1*\edgeDist,0*\edgeDist) circle (2pt) node at (-1*\edgeDist+\putLeft,0*\edgeDist) {$5$};
\filldraw[black] (-1*\edgeDist,-1*\edgeDist) circle (2pt) node at (-1*\edgeDist+\putLeft,-1*\edgeDist) {$6$};

\filldraw[black] (0*\edgeDist,1*\edgeDist) circle (2pt) node at (0,1*\edgeDist+\putLeft) {$7$};
\filldraw[black] (0*\edgeDist,0*\edgeDist) circle (2pt) node at (1*\putLeft,0*\edgeDist) {$8$};

\filldraw[black] (1*\edgeDist,1*\edgeDist) circle (2pt) node at (1*\edgeDist,1*\edgeDist+\putLeft) {$9$};
\filldraw[black] (1*\edgeDist,0*\edgeDist) circle (2pt) node at (1*\edgeDist+\putLeft,0*\edgeDist-\putLeft) {$10$};
\filldraw[black] (1*\edgeDist,-1*\edgeDist) circle (2pt) node at (1*\edgeDist+\putLeft,-1*\edgeDist) {$11$};

\filldraw[black] (2*\edgeDist,1*\edgeDist) circle (2pt) node at (2*\edgeDist,1*\edgeDist+\putLeft) {$12$};
\filldraw[black] (2*\edgeDist,0*\edgeDist) circle (2pt) node at (2*\edgeDist+\putLeft,0*\edgeDist-\putLeft) {$13$};

\filldraw[black] (\iii*\edgeDist,1*\edgeDist) circle (2pt) node at (\iii*\edgeDist,1*\edgeDist+\putLeft) {$14$};
\filldraw[black] (\iii*\edgeDist,0*\edgeDist) circle (2pt) node at (\iii*\edgeDist+\putLeft,0*\edgeDist-\putLeft) {$15$};
\filldraw[black] (\iii*\edgeDist,-1*\edgeDist) circle (2pt) node at (\iii*\edgeDist+\putLeft,-1*\edgeDist) {$16$};

\filldraw[black] (\iv*\edgeDist,1*\edgeDist) circle (2pt) node at (\iv*\edgeDist,1*\edgeDist+\putLeft) {$17$};
\filldraw[black] (\iv*\edgeDist,0*\edgeDist) circle (2pt) node at (\iv*\edgeDist+\putLeft,0*\edgeDist-\putLeft) {$18$};

\filldraw[black] (\v*\edgeDist,1*\edgeDist) circle (2pt) node at (\v*\edgeDist,1*\edgeDist+\putLeft) {$19$};
\filldraw[black] (\v*\edgeDist,0*\edgeDist) circle (2pt) node at (\v*\edgeDist+\putLeft,0*\edgeDist-\putLeft) {$20$};

\filldraw[black] (\vi*\edgeDist,1*\edgeDist) circle (2pt) node at (\vi*\edgeDist,1*\edgeDist+\putLeft) {$21$};
\filldraw[black] (\vi*\edgeDist,0*\edgeDist) circle (2pt) node at (\vi*\edgeDist+\putLeft,0*\edgeDist-\putLeft) {$22$};
\filldraw[black] (\vi*\edgeDist,-1*\edgeDist) circle (2pt) node at (\vi*\edgeDist+\putLeft,-1*\edgeDist) {$23$};

\filldraw[black] (\vii*\edgeDist,1*\edgeDist) circle (2pt) node at (\vii*\edgeDist,1*\edgeDist+\putLeft) {$24$};
\filldraw[black] (\vii*\edgeDist,0*\edgeDist) circle (2pt) node at (\vii*\edgeDist+\putLeft,0*\edgeDist-\putLeft) {$25$};
\filldraw[black] (\vii*\edgeDist,-1*\edgeDist) circle (2pt) node at (\vii*\edgeDist+\putLeft,-1*\edgeDist) {$26$};

\filldraw[black] (\viii*\edgeDist,1*\edgeDist) circle (2pt) node at (\viii*\edgeDist,1*\edgeDist+\putLeft) {$27$};
\filldraw[black] (\viii*\edgeDist,0*\edgeDist) circle (2pt) node at (\viii*\edgeDist+\putLeft,0*\edgeDist-\putLeft) {$28$};
\filldraw[black] (\viii*\edgeDist,-1*\edgeDist) circle (2pt) node at (\viii*\edgeDist+\putLeft,-1*\edgeDist) {$29$};

\draw (-2*\edgeDist,1*\edgeDist) to (-2*\edgeDist,0*\edgeDist); 
\draw (-2*\edgeDist,0*\edgeDist) to (-2*\edgeDist,-1*\edgeDist); 
\draw (-1*\edgeDist,1*\edgeDist) to (-1*\edgeDist,0*\edgeDist); 
\draw (-1*\edgeDist,0*\edgeDist) to (-1*\edgeDist,-1*\edgeDist); 
\draw (0*\edgeDist,1*\edgeDist) to (0*\edgeDist,0*\edgeDist); 
\draw (0*\edgeDist,1*\edgeDist) to (1*\edgeDist,1*\edgeDist); 
\draw (0*\edgeDist,0*\edgeDist) to (1*\edgeDist,1*\edgeDist); 
\draw (1*\edgeDist,1*\edgeDist) to (1*\edgeDist,0*\edgeDist); 
\draw (1*\edgeDist,-1*\edgeDist) to (1*\edgeDist,0*\edgeDist); 
\draw (2*\edgeDist,0*\edgeDist) to (1*\edgeDist,-1*\edgeDist); 
\draw (2*\edgeDist,1*\edgeDist) to (2*\edgeDist,0*\edgeDist); 
\draw (2*\edgeDist,0*\edgeDist) to (3*\edgeDist,0*\edgeDist); 
\draw (2*\edgeDist,1*\edgeDist) to (3*\edgeDist,0*\edgeDist); 
\draw (3*\edgeDist,1*\edgeDist) to (3*\edgeDist,0*\edgeDist); 
\draw (3*\edgeDist,0*\edgeDist) to (4*\edgeDist,0*\edgeDist); 
\draw (3*\edgeDist,0*\edgeDist) to (3*\edgeDist,-1*\edgeDist); 
\draw (4*\edgeDist,1*\edgeDist) to (4*\edgeDist,0*\edgeDist); 
\draw (4*\edgeDist,1*\edgeDist) to (5*\edgeDist,1*\edgeDist); 
\draw (4*\edgeDist,1*\edgeDist) to (5*\edgeDist,0*\edgeDist); 
\draw (4*\edgeDist,0*\edgeDist) to (5*\edgeDist,1*\edgeDist); 
\draw (4*\edgeDist,0*\edgeDist) to (5*\edgeDist,0*\edgeDist); 
\draw (5*\edgeDist,1*\edgeDist) to (5*\edgeDist,0*\edgeDist); 
\draw (5*\edgeDist,1*\edgeDist) to (6*\edgeDist,1*\edgeDist); 
\draw (6*\edgeDist,0*\edgeDist) to (6*\edgeDist,-1*\edgeDist); 
\draw (6*\edgeDist,1*\edgeDist) to (6*\edgeDist,0*\edgeDist); 
\draw (6*\edgeDist,1*\edgeDist) to (7*\edgeDist,0*\edgeDist); 
\draw (6*\edgeDist,0*\edgeDist) to (7*\edgeDist,0*\edgeDist); 
\draw (7*\edgeDist,1*\edgeDist) to (7*\edgeDist,0*\edgeDist); 
\draw (7*\edgeDist,0*\edgeDist) to (8*\edgeDist,0*\edgeDist); 
\draw (7*\edgeDist,0*\edgeDist) to (8*\edgeDist,1*\edgeDist); 
\draw (7*\edgeDist,0*\edgeDist) to (7*\edgeDist,-1*\edgeDist); 
\draw (8*\edgeDist,1*\edgeDist) to (8*\edgeDist,0*\edgeDist); 
\draw (8*\edgeDist,0*\edgeDist) to (8*\edgeDist,-1*\edgeDist); 

\draw (-2*\edgeDist,0*\edgeDist) to (-1*\edgeDist,0*\edgeDist); 
\draw (-2*\edgeDist,-1*\edgeDist) to (-1*\edgeDist,0*\edgeDist); 
\draw (-1*\edgeDist,0*\edgeDist) to[out=-45,in=225] (1*\edgeDist,0*\edgeDist);
\end{tikzpicture}
\end{center}
\caption{Block Graph}
\label{fig_forbidden_paths}
\end{figure}

\subsection{DOIP Property for Paths and Forbidden Subgraphs}

The main result of this subsection is the following theorem, which gives a combinatorial characterization when $\ul{P} = P_{1},\ldots,P_{\ell}$, vertex-disjoint induced paths of a block graph, are DOIP. This characterization does not refer to the labeling of the paths, $\sigma$, or to the orientations of the paths $P_{i}$. In the subsequent section, we will leverage this theorem to prove the equality of $\nu(G)$ and $\reg(S/J_{G})$ for block graphs.

\bthm
\label{thm_doip_iff_no_forbidden_path}
Let $G$ be a block graph, and let $\ul{P} := P_{1},\ldots,P_{\ell}$ be vertex-disjoint induced paths of $G$. The following statements are equivalent:
\begin{enumerate}
    \item $\ul{P}$ is DOIP for any choice of orientations $\phi_{i}$, $1 \leq i \leq \ell$,
    \item $P_{\ind}$ does not contain an internal strand, an internal fork, a double fork, or a complete ladder as subgraphs.
\end{enumerate}
\ethm

The idea for the proof of Theorem \ref{thm_doip_iff_no_forbidden_path} is as follows: For (1) implies (2), we show that if $P_{\ind}$ contains an internal strand, an internal fork, a double fork, or a complete ladder as a subgraph, then $K_{P_{\ind}}$ contains a directed two cycle. For (2) implies (1), it suffices to show that $K_{P_{\ind}}$ is directed acyclic. We would like to say that if $K_{P_{\ind}}$ has a directed cycle, then the paths realizing this directed cycle of $K_{P_{\ind}}$ would realize a \textit{large} cycle in $G$. This would lead to complications, since the induced subgraph on the vertices of any cycle in a block graph is a complete graph. The difficulty is that the paths realizing the directed cycle of $K_{P_{\ind}}$ may a priori intersect each other, thwarting this hope. However, our observation is that these paths will intersect at the cost of introducing an internal strand, an internal fork, a double fork, or a complete ladder as a subgraph in $G$.

We begin by showing that $K_{P_{\ind}}$ is a simple multigraph  for any block graph.

\bprop 
\label{prop_no_1_cycle}
Let $\ul{P} := P_{1},\ldots,P_{\ell}$ be vertex-disjoint induced oriented paths of a block graph $G$. Then, $K_{P_{\ind}}$ has no loops or multiarcs, i.e., $K_{P_{\ind}}$ is a (simple) directed graph.
\eprop

\bproof
Suppose by contradiction that $K_{P_{\ind}}$ has a loop or a multiarc. If $K_{P_{\ind}}$ has a loop, then there is an oriented induced path $Q$ of $G$ such that $\phi_{Q}(1) = \phi_{i}(1)$ and $\phi_{Q}(2) = \phi_{i}(2)$ for some $1 \leq i \leq \ell$. If $K_{P_{\ind}}$ has a multiarc, then there are oriented induced paths $Q_{1}$ and $Q_{2}$ of $G$ such that $\phi_{Q_{1}}(j) = \phi_{Q_{2}}(j)$ for $j = 1,2$. Lemma \ref{lem_decomp_intersecting_induced_paths} implies that the intersection of any two induced paths of a block graph consists of exactly one connected component. The only way that two induced paths can have the same terminal vertices is if they are the same path. Thus, $Q = P_{i}$ and $Q_{1} = Q_{2}$, a contradiction.
\eproof

The following proposition connects arcs of $K_{P_{\ind}}$ to strands of $P_{\ind}$.

\bprop
\label{prop_arc_iff_strand}
$K_{P_{\ind}}$ contains the arc $(i,j)$ if and only if $P_{\ind}$ contains a strand from $P_{i}$ to $P_{j}$.
\eprop

\bproof
$(\implies)$ There exists an induced path $Q$ from $\phi_{i}(1)$ to $\phi_{j}(2)$ realizing the arc $(i,j)$ of $K_{P_{\ind}}$. Lemma \ref{lem_decomp_induced_path_intersecting_two_others} applied to $Q$, $P_{i}$, and $P_{j}$ implies that $Q = Q_{1} \ast \gamma \ast Q_{2}$. Then, $\gamma$ is a strand from $P_{i}$ to $P_{j}$.

$(\impliedby)$ Let $Q$ be a strand from $P_{i}$ to $P_{j}$. Let $R_{i}$ (respectively, $R_{j}$) be the subpath of $P_{i}$ (respectively, from $P_{j}$) from $\phi_{i}(1)$ to $\phi_{Q}(1)$ (respectively, $\phi_{Q}(2)$ to $\phi_{j}(2)$). Let $T$ be the induced path from $\phi_{i}(1)$ to $\phi_{j}(2)$. In order to show that $T$ realizes the arc $(i,j)$ of $K_{P_{\ind}}$, it suffices to show that $T$ does not contain $P_{k}$ for $1 \leq k \leq \ell$. We observe that since $R_{i}\ast Q \ast R_{j}$ is a path from $\phi_{i}(1)$ to $\phi_{j}(2)$, we have that
\begin{align*}
    V(T) 
    &\subset R_{i}\ast Q \ast R_{j} \\ 
    &\subset \left( V(P_{i}) \setminus \{\phi_{i}(2)\} \right) \cup V(Q) \cup \left( V(P_{j}) \setminus \{\phi_{j}(1)\} \right).
\end{align*}
We observe that
\begin{enumerate}
    \item $V(P_{i}) \not\subset V(T)$ because $\phi_{i}(2) \notin V(Q)$, as $Q$ is a strand,
    \item $V(P_{j}) \not\subset V(T)$ because $\phi_{j}(1) \notin V(Q)$, as $Q$ is a strand,
    \item $V(P_{k}) \not\subset V(T)$ for $k \in [\ell] \setminus \{i,j\}$; otherwise, $V(P_{k}) \subset V(Q)$ (since the $P_{i}$ are vertex-disjoint), which is impossible because $Q$ is a strand.
\end{enumerate}
\eproof

The next result establishes that (1) implies (2) of Theorem \ref{thm_doip_iff_no_forbidden_path}.

\bcor 
\label{cor_KP_ind_2_cycle_iff_Pind_weakly_forbidden_path}
Let $\ul{P} := P_{1},\ldots,P_{\ell}$ be vertex-disjoint oriented induced paths of a block graph $G$. Then, $P_{\ind}$ contains an internal strand, an internal fork, a double fork, or a complete ladder if and only if $K_{P_{\ind}}$ has a directed cycle of length two.
\ecor

\bproof
$(\implies)$ Proposition \ref{prop_arc_iff_strand} implies that it suffices to construct a strand from $P_{i}$ to $P_{j}$ and a strand from $P_{j}$ to $P_{i}$ whenever $P_{\ind}$ contains one of the graphs in question. 

If $Q$ is an internal strand from $P_{i}$ to $P_{j}$, then it is clear that $Q$ is a strand from $P_{i}$ to $P_{j}$ and vice versa. 

Let $H$ be an internal fork from $P_{i}$ to $P_{j}$. Without loss of generality, we may suppose that $b \neq \phi_{j}(2)$ and that $c \neq \phi_{j}(1)$. Then, $Q \ast [\phi_{Q}(2),c]$ is a strand from $P_{i}$ to $P_{j}$, and $[b,\phi_{Q}(2)] \ast Q$ is a strand from $P_{j}$ to $P_{i}$. 

Let $H$ be a double fork from $P_{i}$ to $P_{j}$. Without loss of generality, we may suppose that $a \neq \phi_{i}(2)$, $b \neq \phi_{i}(1)$, $c \neq \phi_{j}(2)$, and $d \neq \phi_{j}(1)$. Then, $[a,\phi_{Q}(1)]\ast Q \ast [\phi_{Q}(2),d]$ is a strand from $P_{i}$ to $P_{j}$, and $[c,\phi_{Q}(2)] \ast Q \ast [\phi_{Q}(1),b]$ is a strand from $P_{j}$ to $P_{i}$. 

Suppose that $H$ is a complete ladder. Without loss of generality, we may suppose that $a$ (respectively, $c$) is closer to $\phi_{i}(1)$ (respectively, $\phi_{j}(1)$) than $b$ (respectively, $d$). Then, $[a,d]$ is a strand from $P_{i}$ to $P_{j}$, and $[c,b]$ is a strand from $P_{j}$ to $P_{i}$.

$(\impliedby$) Suppose that $K_{P_{\ind}}$ has a directed cycle of length two; then there exist induced paths $Q_{1}$ and $Q_{2}$ of $P_{\ind}$ realizing this directed cycle. Without loss of generality, we may suppose that $Q_{1}$ (respectively, $Q_{2}$) realizes the arc $(1,2)$ (respectively, $(2,1)$). Lemma \ref{lem_decomp_induced_path_intersecting_two_others} implies that $Q_{i} = Q_{i,1} \ast \gamma_{i} \ast Q_{i,2}$ for $i = 1,2$ where $Q_{i,1}$ is a subpath of $P_{1}$ and $Q_{i,2}$ is a subpath of $P_{2}$. Let $a$ (respectively $b$) be the terminal vertex of $\gamma_{1}$ (respectively $\gamma_{2}$) contained in $P_{1}$. Let $c$ (respectively $d$) be the terminal vertex of $\gamma_{1}$ (respectively $\gamma_{2}$) contained in $P_{2}$. We denote the subpath of $P_{1}$ (respectively $P_{2}$) having terminal vertices $a$ and $b$ (respectively $c$ and $d$) by $[a,b]$ (respectively $[c,d]$).

\ul{Case 1.} Suppose that $V(\gamma_{1}) \cap V(\gamma_{2}) = \varnothing$. Then, $a$, $b$, $c$, and $d$ are distinct vertices, and
\begin{align*}
    [a,b] \ast \gamma_{2} \ast [d,c] \ast \gamma_{1}
\end{align*}
is a cycle. Since $G$ is a block graph, the induced subgraph on $\{a,b,c,d\}$ is a complete graph. Thus, $P_{\ind}$ contains a complete ladder.

\ul{Case 2.} Refer to Figure \ref{fig_KP_ind_2_cycle_iff_Pind_weakly_forbidden_path} for this case. Suppose that $V(\gamma_{1}) \cap V(\gamma_{2}) \neq \varnothing$. Lemma \ref{lem_decomp_induced_path_intersecting_two_others} applied to $\gamma_{1}$ and $\gamma_{2}$ implies that there is a decomposition $\gamma_{i} = \mu_{i,1} \ast \omega \ast \mu_{i,2}$ where $\mu_{i,1}$ intersects $P_{1}$ and $\mu_{i,2}$ intersects $P_{2}$ for $i = 1,2$.

When $a = b$ and $c = d$, $\omega$ is an internal path. When $a = b$ and $c \neq d$, the subgraph of $P_{\ind}$ having vertices $V(\omega) \cup \{c,d\}$ and edges $E(\omega) \cup E(\mu_{1,2}) \cup E(\mu_{2,2})$ is an internal fork. When $a \neq b$ and $c \neq d$, the subgraph of $P_{\ind}$ having vertices $V(\omega) \cup \{a,b,c,d\}$ and edges $E(\omega) \cup E(\mu_{1,1}) \cup E(\mu_{1,2}) \cup E(\mu_{2,1}) \cup E(\mu_{2,2})$ is a double fork.
\eproof

\def\putLeft{.5*.8}
\def\edgeDist{1.75*.8}
\begin{figure}
\begin{center}
\begin{tikzpicture}
\filldraw[black] (0*\edgeDist,0*\edgeDist) circle (2pt);
\filldraw[black] (0*\edgeDist,1*\edgeDist) circle (2pt) node at (-\putLeft,1*\edgeDist) {$b$};
\filldraw[black] (0*\edgeDist,2*\edgeDist) circle (2pt) node at (-\putLeft,2*\edgeDist) {$a$};
\filldraw[black] (0*\edgeDist,3*\edgeDist) circle (2pt);
\filldraw[black] (0*\edgeDist,2*\edgeDist) circle (0pt) node at (-\putLeft,.5*\edgeDist) {$\vdots$};
\filldraw[black] (0*\edgeDist,2*\edgeDist) circle (0pt) node at (-\putLeft,2.7*\edgeDist) {$\vdots$};
\filldraw[black] (1*\edgeDist,1.5*\edgeDist) circle (2pt) node at (1.2*\edgeDist,1.5*\edgeDist-\putLeft) {$\omega_{1}$};
\filldraw[black] (4*\edgeDist,1.5*\edgeDist) circle (2pt) node at (4*\edgeDist,1.5*\edgeDist-\putLeft) {$\omega_{r-1}$};
\filldraw[black] (2.5*\edgeDist,1.5*\edgeDist) circle (0pt) node at (2.5*\edgeDist,1.5*\edgeDist+\putLeft) {$\omega$};
\filldraw[black] (2.5*\edgeDist,1.5*\edgeDist) circle (0pt) node at (2.5*\edgeDist,1.5*\edgeDist-\putLeft) {$\cdots$};
\filldraw[black] (5*\edgeDist,0*\edgeDist) circle (2pt);
\filldraw[black] (5*\edgeDist,1*\edgeDist) circle (2pt) node at (5*\edgeDist+\putLeft,1*\edgeDist) {$c$};
\filldraw[black] (5*\edgeDist,2*\edgeDist) circle (2pt) node at (5*\edgeDist+\putLeft,2*\edgeDist) {$d$};
\filldraw[black] (5*\edgeDist,3*\edgeDist) circle (2pt);
\filldraw[black] (5*\edgeDist,2*\edgeDist) circle (0pt) node at (5*\edgeDist+\putLeft,.5*\edgeDist) {$\vdots$};
\filldraw[black] (5*\edgeDist,2*\edgeDist) circle (0pt) node at (5*\edgeDist+\putLeft,2.7*\edgeDist) {$\vdots$};
\filldraw[black] (.5*\edgeDist,1.75*\edgeDist) circle (0pt) node at (.5*\edgeDist+.5*\putLeft,1.6*\edgeDist+\putLeft) {$\mu_{1,1}$};
\filldraw[black] (.5*\edgeDist,1.75*\edgeDist) circle (0pt) node at (.5*\edgeDist+.5*\putLeft,1.3*\edgeDist-\putLeft) {$\mu_{2,1}$};
\filldraw[black] (4.5*\edgeDist,1.75*\edgeDist) circle (0pt) node at (4.5*\edgeDist+.5*\putLeft,1.75*\edgeDist+\putLeft) {$\mu_{2,2}$};
\filldraw[black] (4.5*\edgeDist,1.75*\edgeDist) circle (0pt) node at (4.5*\edgeDist+.5*\putLeft,1.2*\edgeDist-\putLeft) {$\mu_{1,2}$};
%
%
%
\draw (0*\edgeDist,0*\edgeDist) to (0*\edgeDist,1*\edgeDist); 
\draw (0*\edgeDist,1*\edgeDist) to (0*\edgeDist,2*\edgeDist); 
\draw (0*\edgeDist,2*\edgeDist) to (0*\edgeDist,3*\edgeDist); 

\draw (5*\edgeDist,0*\edgeDist) to (5*\edgeDist,1*\edgeDist); 
\draw (5*\edgeDist,1*\edgeDist) to (5*\edgeDist,2*\edgeDist); 
\draw (5*\edgeDist,2*\edgeDist) to (5*\edgeDist,3*\edgeDist); 

\draw (0*\edgeDist,2*\edgeDist) to (1*\edgeDist,1.5*\edgeDist); 
\draw (0*\edgeDist,1*\edgeDist) to (1*\edgeDist,1.5*\edgeDist); 
\draw (5*\edgeDist,2*\edgeDist) to (4*\edgeDist,1.5*\edgeDist); 
\draw (5*\edgeDist,1*\edgeDist) to (4*\edgeDist,1.5*\edgeDist); 
\draw (1*\edgeDist,1.5*\edgeDist) to (4*\edgeDist,1.5*\edgeDist); 
\end{tikzpicture}
\end{center}
\caption{Illustration for Corollary \ref{cor_KP_ind_2_cycle_iff_Pind_weakly_forbidden_path}}
\label{fig_KP_ind_2_cycle_iff_Pind_weakly_forbidden_path}
\end{figure}

The following lemmas will help us control the intersection of strands.

\blem 
\label{lem_strand_intersecting_path_creates_new_arc}
If $\gamma$ is a strand from $P_{i}$ to $P_{k}$ and $V(\gamma) \cap V(P_{j}) \neq \varnothing$ for some $j \in [\ell] \setminus \{i,k\}$, then $K_{P_{\ind}}$ contains the arc $(i,j)$ or the arc $(j,k)$.
\elem 

\bproof
By Lemma \ref{lem_decomp_intersecting_induced_paths}, we can write $\gamma$ as $\gamma_{1} \ast \omega \ast \gamma_{2}$ where $V(\gamma_{n}) \cap V(P_{j}) = \{v_{n}\}$ for some vertices $v_{n}$ for $n = 1,2$, $\omega$ is a subpath of $P_{j}$, $V(\gamma_{1}) \cap V(P_{i}) \neq \varnothing$, and $V(\gamma_{2}) \cap V(P_{k}) \neq \varnothing$. If $v_{1} \neq \phi_{j}(1)$, then $\gamma_{1}$ is a strand from $P_{i}$ to $P_{j}$. If $v_{2} \neq \phi_{j}(2)$, then $\gamma_{2}$ is a strand from $P_{j}$ to $P_{k}$. It cannot be the case that both $v_{1} = \phi_{j}(1)$ and $v_{2} = \phi_{j}(2)$, as $\gamma$ does not contain $P_{j}$.
\eproof

\bprop 
\label{prop_intersection_two_strands_is_strand}
Suppose that $P_{\ind}$ does not contain an internal strand, an internal fork, a double fork, or a complete ladder.
Let $\gamma_{1}$ be a strand from $P_{a}$ to $P_{b}$, and $\gamma_{2}$ be a strand from $P_{b}$ to $P_{c}$. If $V(\gamma_{1}) \cap V(\gamma_{2}) \neq \varnothing$, then $K_{P_{\ind}}$ contains the arc $(a,c)$.
\eprop 

\bproof
For this proof, refer to Figure \ref{fig_intersection_two_strands_is_strand}. 
If $V(\gamma_{1}) \cap V(P_{c}) \neq \varnothing$, then Lemma \ref{lem_strand_intersecting_path_creates_new_arc} implies that $(a,c)$ or $(c,b)$ is an arc of $K_{P_{\ind}}$. Corollary \ref{cor_KP_ind_2_cycle_iff_Pind_weakly_forbidden_path} implies that $(c,b)$ cannot be an arc of $K_{P_{\ind}}$. Likewise, it is shown that if $V(\gamma_{2}) \cap V(P_{a}) \neq \varnothing$, then $(a,c)$ is an arc of $K_{P_{\ind}}$. Thus, we may assume that $V(\gamma_{1}) \cap V(P_{c}) = \varnothing$ and that  $V(\gamma_{2}) \cap V(P_{a}) = \varnothing$. Lemma \ref{lem_decomp_induced_path_intersecting_two_others} implies that there is a decomposition 
\begin{align*}
    \gamma_{i} = \mu_{i,1} \ast \omega \ast \mu_{i,2},
\end{align*}
where $\mu_{i,j}$ is determined by the condition of containing the vertex $\phi_{\gamma_{i}}(j)$ for $1 \leq i,j \leq 2$. We denote the intersection of $\mu_{2,1}$ (respectively, $\mu_{1,2}$) with $P_{b}$ by the vertex $t$ (respectively, $s$). We denote the vertex that is the intersection of $\mu_{2,1}$ and $\mu_{1,2}$ by $r$, and we denote the vertex which is the intersection of $\mu_{1,1}$ and $\mu_{2,2}$ by $v$. We denote the vertices belonging to $\mu_{1,1}$ and $\mu_{2,2}$ and adjacent to $v$ by $u$ and $w$, respectively. We observe that either $r = s = t$ (in which case $r$ is an internal vertex of $P_{b}$), or the induced subgraph on $\{r,s,t\}$ is a complete graph (by Lemma \ref{lem_block_graphs_do_not_have_certain_cycles}). We consider the path $Q = \mu_{1,1} \ast \mu_{2,2}$. We observe that $Q$ is not an induced subpath of $P_{\ind}$ if and only if $\{u,w\} \in E(P_{\ind})$ (by Lemma \ref{lem_block_graphs_do_not_have_certain_cycles}). Let $\tilde{Q}$ denote the induced subpath of $Q$. We show that $\tilde{Q}$ does not contain $P_{k}$ for any $1\leq k \leq \ell$. Otherwise, it would follow that $\tilde{Q}$ is a strand from $P_{a}$ to $P_{c}$, and the result would follow from Proposition \ref{prop_arc_iff_strand}. 

Suppose by contradiction that $\tilde{Q}$ contains $P_{k}$ for some $1 \leq k \leq \ell$. Since $\gamma_{1}$ does not intersect $P_{c}$ and $\gamma_{2}$ does not intersect $P_{a}$, it follows that $\tilde{Q}$ does not contain $P_{a}$ or $P_{c}$.

\ul{Case 1.} Suppose that $\tilde{Q} = Q$ is an induced path. Then, $v$ is an internal vertex of $P_{k}$; otherwise, $P_{k}$ would be properly contained in either $\mu_{1,1}$ or in $\mu_{2,2}$, which would contradict $\gamma_{1}$ and $\gamma_{2}$ being strands. When $s \neq t$, the subgraph having vertices $V(\omega) \cup \{s,t\}$ and edges $E(\omega) \cup \{ \{r,s\}, \{r,t\} \}$ is an internal fork. When $s = t$, the subgraph having vertices $V(\omega) \cup \{s\}$ and edges $E(\omega) \cup \{ \{r,s\} \}$ is an internal strand.

\ul{Case 2.} Suppose that $\tilde{Q} \neq Q$. Then, $\tilde{Q}$ contains the edge $\{u,w\}$. Moreover, $P_{k}$ contains the edge $\{u,w\}$, since $P_{k}$ is not contained in $\gamma_{1}$ or $\gamma_{2}$. It follows that there is a subgraph $H$ with $V(H) \subset V(\omega) \cup \{u,w\} \cup \{s,t\}$ and with $E(H) \subset E(\omega) \cup \{ \{u,v\},\{r,s\},\{r,t\} \}$, which is an internal strand, an internal fork, or a double fork.
\eproof

\def\putLeft{.5*.8}
\def\edgeDist{1.75*.8}
\begin{figure}
\begin{center}
\begin{tikzpicture}[decoration={brace,amplitude=15pt}]
\filldraw[black] (0*\edgeDist,0*\edgeDist) circle (2pt) node at (-\putLeft,0*\edgeDist-\putLeft) {$P_{c}$};
\filldraw[black] (0*\edgeDist,1*\edgeDist) circle (2pt) node at (-\putLeft,1*\edgeDist) {$\vdots$};
\filldraw[black] (0*\edgeDist,2*\edgeDist) circle (2pt);
\filldraw[black] (0*\edgeDist,3*\edgeDist) circle (2pt);
\filldraw[black] (0*\edgeDist,4*\edgeDist) circle (2pt) node at (-\putLeft,4*\edgeDist) {$\vdots$};
\filldraw[black] (0*\edgeDist,5*\edgeDist) circle (2pt) node at (0*\edgeDist-\putLeft,5*\edgeDist+\putLeft) {$P_{a}$};
\filldraw[black] (2*\edgeDist,2*\edgeDist) circle (2pt) node at (2*\edgeDist-\putLeft,2*\edgeDist) {$w$};
\filldraw[black] (2*\edgeDist,3*\edgeDist) circle (2pt) node at (2*\edgeDist-\putLeft,3*\edgeDist) {$u$};
\filldraw[black] (3*\edgeDist,2.5*\edgeDist) circle (2pt) node at (3*\edgeDist,2.5*\edgeDist-\putLeft) {$v$};
\filldraw[black] (2*\edgeDist,2*\edgeDist) circle (0pt) node at (1.5*\edgeDist+\putLeft,1.75*\edgeDist-2*\putLeft) {$\mu_{2,2}$};
\filldraw[black] (2*\edgeDist,2*\edgeDist) circle (0pt) node at (1.5*\edgeDist+\putLeft,2.75*\edgeDist+3.5*\putLeft) {$\mu_{1,1}$};
\filldraw[black] (2*\edgeDist,2*\edgeDist) circle (0pt) node at (5.5*\edgeDist,3*\edgeDist) {$\mu_{2,1}$};
\filldraw[black] (2*\edgeDist,2*\edgeDist) circle (0pt) node at (5.5*\edgeDist,2*\edgeDist) {$\mu_{1,2}$};
\filldraw[black] (2*\edgeDist,2*\edgeDist) circle (0pt) node at (4*\edgeDist,2.5*\edgeDist+\putLeft) {$\omega$};
\filldraw[black] (5*\edgeDist,2.5*\edgeDist) circle (2pt) node at (5*\edgeDist,2.5*\edgeDist+\putLeft) {$r$};
\filldraw[black] (6*\edgeDist,1*\edgeDist) circle (2pt) node at (6*\edgeDist+\putLeft,1*\edgeDist-\putLeft) {$P_{b}$};
\filldraw[black] (6*\edgeDist,2*\edgeDist) circle (2pt) node at (6*\edgeDist+\putLeft,2*\edgeDist) {$s$};
\filldraw[black] (6*\edgeDist,3*\edgeDist) circle (2pt) node at (6*\edgeDist+\putLeft,3*\edgeDist) {$t$};
\filldraw[black] (6*\edgeDist,4*\edgeDist) circle (2pt);
\filldraw[black] (6*\edgeDist+\putLeft,1.5*\edgeDist) circle (0pt) node at (6*\edgeDist+\putLeft,1.5*\edgeDist) {$\vdots$};
\filldraw[black] (6*\edgeDist+\putLeft,1.5*\edgeDist) circle (0pt) node at (6*\edgeDist+\putLeft,3.5*\edgeDist) {$\vdots$};
%
%
%
\draw (0*\edgeDist,0*\edgeDist) to (0*\edgeDist,1*\edgeDist); 
\draw (0*\edgeDist,1*\edgeDist) to (0*\edgeDist,2*\edgeDist); 
\draw (0*\edgeDist,3*\edgeDist) to (0*\edgeDist,4*\edgeDist); 
\draw (0*\edgeDist,4*\edgeDist) to (0*\edgeDist,5*\edgeDist); 
\draw (6*\edgeDist,1*\edgeDist) to (6*\edgeDist,2*\edgeDist); 
\draw (6*\edgeDist,2*\edgeDist) to (6*\edgeDist,3*\edgeDist); 
\draw (6*\edgeDist,3*\edgeDist) to (6*\edgeDist,4*\edgeDist); 
\draw (0*\edgeDist,1*\edgeDist) to (3*\edgeDist,2.5*\edgeDist); 
\draw (0*\edgeDist,4*\edgeDist) to (3*\edgeDist,2.5*\edgeDist); 
\draw (3*\edgeDist,2.5*\edgeDist) to (5*\edgeDist,2.5*\edgeDist); 
\draw (5*\edgeDist,2.5*\edgeDist) to (6*\edgeDist,2*\edgeDist); 
\draw (5*\edgeDist,2.5*\edgeDist) to (6*\edgeDist,3*\edgeDist); 
\draw[dashed] (2*\edgeDist,2*\edgeDist) to (2*\edgeDist,3*\edgeDist); 
\draw[decorate]  (3*\edgeDist,2.5*\edgeDist) -- (0*\edgeDist,1*\edgeDist); 
\draw[decorate]  (0*\edgeDist,4*\edgeDist) -- (3*\edgeDist,2.5*\edgeDist); 
\end{tikzpicture}
\end{center}
\caption{Illustration for Proposition \ref{prop_intersection_two_strands_is_strand}}
\label{fig_intersection_two_strands_is_strand}
\end{figure}

\bprop 
\label{prop_intersection_two_strands_and_four_paths_is_strand}
Suppose that $P_{\ind}$ does not contain an internal strand, an internal fork, a double fork, or a complete ladder.
Let $\gamma_{1}$ be a strand from $P_{a}$ to $P_{b}$, and $\gamma_{2}$ be a strand from $P_{c}$ to $P_{d}$ where $a,b,c,d \in [\ell]$ are distinct. If $V(\gamma_{1}) \cap V(\gamma_{2}) \neq \varnothing$, then $K_{P_{\ind}}$ contains the arc $(a,c)$, $(a,d)$, $(b,d)$, $(c,a)$, $(c,b)$, $(d,a)$, or $(d,b)$.

If, in addition $V(\gamma_{1}) \cap V(P_{i}) = \varnothing$ for $i \in \{c,d\}$ and $V(\gamma_{2}) \cap V(P_{j}) = \varnothing$ for $j \in \{a,b\}$, then $K_{P_{\ind}}$ contains the arc $(a,d)$ or the arc $(c,b)$.
\eprop 

\bproof
For this proof, refer to Figure \ref{fig_intersection_two_strands_and_four_paths_is_strand}. By Proposition \ref{lem_strand_intersecting_path_creates_new_arc}, we may assume that $V(\gamma_{1}) \cap P_{i} = \varnothing$ for $i \in \{c,d\}$ and $V(\gamma_{2}) \cap P_{i} = \varnothing$ for $i \in \{a,b\}$. Lemma \ref{lem_decomp_induced_path_intersecting_two_others} implies that there is a decomposition 
\begin{align*}
    \gamma_{i} = \mu_{i,1} \ast \omega \ast \mu_{i,2},
\end{align*}
where $\mu_{i,j}$ is determined by the condition of containing the vertex $\phi_{\gamma_{i}}(j)$ for $1 \leq i,j \leq 2$.

\ul{Case 1.} Suppose that $V(\mu_{1,1}) \cap V(\mu_{2,2}) \neq \varnothing$. Then, we have that $V(\mu_{1,1}) \cap V(\mu_{2,2}) = \{v\}$ and that $V(\mu_{1,2}) \cap V(\mu_{2,1}) = r$. (When $\card{V(\omega)} = 1$, we have that $v = r$.) We denote by $u$ and $t$ the vertices of $\mu_{1,1}$ and $\mu_{2,2}$, respectively, which are adjacent to $v$. We denote by $s$ and $w$ the vertices of $\mu_{1,2}$ and $\mu_{2,1}$, respectively, which are adjacent to $r$. We denote by $Q_{1}$ and $Q_{2}$ the paths $\mu_{1,1}\ast\mu_{2,2}$ and $\mu_{2,1}\ast\mu_{1,2}$, respectively. The paths $Q_{1}$ and $Q_{2}$ are not an induced path if and only if $\{u,t\}$ and $\{w,s\}$ are induced edges of $Q_{1}$ and $Q_{2}$, respectively. Let $\tilde{Q}_{1}$ and $\tilde{Q}_{2}$ denote the induced subpath of $P_{\ind}$ on $V(Q_{1})$ and $V(Q_{2})$, respectively. We will show that it is not possible for both $\tilde{Q}_{1}$ and $\tilde{Q}_{2}$ to contain paths $P_{i}$ and $P_{j}$ for some $i,j \in [\ell]$ distinct. In which case, it follows that $\tilde{Q_{1}}$ or $\tilde{Q}_{2}$ is a strand from $P_{a}$ to $P_{d}$ or a strand from $P_{c}$ to $P_{b}$, respectively. The result then follows from Proposition \ref{prop_arc_iff_strand}. 

Suppose by contradiction that $\tilde{Q}_{1}$ contains $P_{i}$ and that $\tilde{Q}_{2}$ contains $P_{j}$, and we consider the following subcases.

\ul{Subcase 1.(a).} Suppose that $\tilde{Q}_{1} = Q_{1}$ and that $\tilde{Q}_{2} = Q_{2}$. First, we observe that $P_{i}$ contains $v$ as an internal vertex; otherwise, $P_{i}$ would belong to $\gamma_{1}$ or $\gamma_{2}$. For similar reasons, $P_{j}$ contains $r$ as an internal vertex. When $\card{V(\omega)} = 1$, $r = v$; contradicting $P_{i}$ and $P_{j}$ being vertex-disjoint. When $V(\omega) \geq 2$, $\omega$ would be an internal strand, a contradiction.

\ul{Subcase 1.(b).} Suppose that $\tilde{Q_{1}}$ contains the edge $\{u,t\}$ and that $\tilde{Q}_{2} = Q_{2}$. Then, $\{u,t\} \in E(P_{i})$, and $r$ is an internal vertex of $P_{j}$. It follows that the subgraph having vertices $V(\omega) \cup \{u,t\}$ and edges $E(\omega) \cup \{ \{u,v\},\{t,v\} \}$ is an internal fork, a contradiction.

\ul{Subcase 1.(c).} Suppose that $\tilde{Q}_{1} = Q_{1}$ and that $\tilde{Q}_{2}$ contains the edge $\{w,s\}$. This subcase is analogous to subcase 1.(b).

\ul{Subcase 1.(d).} Suppose that $\tilde{Q_{1}}$ contains the edge $\{u,t\}$ and that $\tilde{Q_{2}}$ contains the edge $\{w,s\}$. Then, $\{u,t\} \in E(P_{i})$ and $\{w,s\} \in E(P_{j})$. Consequently, the subgraph having vertices $V(\omega) \cup \{u,t,w,s\}$ and edges $E(\omega) \cup \{ \{u,v\}, \{t,v\}, \{w,r\}, \{s,r\} \}$ is a double fork, a contradiction.

\ul{Case 2.} Suppose that $V(\mu_{1,1}) \cap V(\mu_{2,1}) \neq \varnothing$. Then, we have that $V(\mu_{1,1}) \cap V(\mu_{2,1}) = \{v\}$ and that $V(\mu_{1,2}) \cap V(\mu_{2,2}) = \{r\}$. (We may suppose that $v \neq r$; otherwise, we would be in Case 1.) We denote by $u$ and $w$ the vertices of $\mu_{1,1}$ and $\mu_{2,1}$, respectively, which are adjacent to $v$. We denote by $s$ and $t$ the vertices of $\mu_{1,2}$ and $\mu_{2,2}$, respectively, which are adjacent to $r$. We denote by $Q_{1}$ and $Q_{2}$ the paths $\mu_{1,1}\ast\omega\ast\mu_{2,2}$ and $\mu_{2,1}\ast\omega\ast\mu_{1,2}$, respectively. We observe that $Q_{1}$ is an induced path of $P_{\ind}$ by Lemma \ref{lem_block_graphs_do_not_have_certain_cycles} together with the observations that $\mu_{1,1}\ast\omega$ and $\omega\ast\mu_{2,2}$ are induced paths of $P_{\ind}$ being subpaths of the induced paths $\gamma_{1}$ and $\gamma_{2}$, respectively. Similarly, $Q_{2}$ is an induced path of $P_{\ind}$. If $Q_{1}$ contains the path $P_{i}$ for some $i \in [\ell]$, then $u$ and $t$ belong to $V(P_{i})$; otherwise, $P_{i}$ would be contained in $\gamma_{1}$ or $\gamma_{2}$. It follows that $Q_{2}$ cannot contain any path $P_{j}$ for $j \in [\ell]$ as such a path would necessarily be vertex-disjoint from $P_{i}$ and contain $\omega$, which is impossible.
\eproof 

\def\putLeft{.5*.8}
\def\edgeDist{1.75*.8}
\def\horDist{1.5}
\begin{figure}
\begin{subfigure}{.45\linewidth}
\begin{tikzpicture}[decoration={brace,amplitude=22pt}]
\filldraw[black] (0*\edgeDist,0*\edgeDist) circle (2pt) node at (-\putLeft,0*\edgeDist-\putLeft) {$P_{c}$};
\filldraw[black] (0*\edgeDist,1*\edgeDist) circle (2pt) node at (-\putLeft,1*\edgeDist) {$\vdots$};
\filldraw[black] (0*\edgeDist,2*\edgeDist) circle (2pt);
\filldraw[black] (0*\edgeDist,3*\edgeDist) circle (2pt);
\filldraw[black] (0*\edgeDist,4*\edgeDist) circle (2pt) node at (-\putLeft,4*\edgeDist) {$\vdots$};
\filldraw[black] (0*\edgeDist,5*\edgeDist) circle (2pt) node at (0*\edgeDist-\putLeft,5*\edgeDist+\putLeft) {$P_{a}$};
\filldraw[black] (3.5*\edgeDist,2*\edgeDist) circle (2pt) node at (3.5*\edgeDist+.75*\putLeft,2*\edgeDist) {$s$};
\filldraw[black] (2*\edgeDist,2*\edgeDist) circle (2pt) node at (2*\edgeDist,2*\edgeDist-.6*\putLeft) {$w$};
\filldraw[black] (2.5*\edgeDist,3*\edgeDist) circle (2pt) node at (2.5*\edgeDist,3*\edgeDist+.6*\putLeft) {$t$};
\filldraw[black] (1*\edgeDist,3*\edgeDist) circle (2pt) node at (1*\edgeDist,3*\edgeDist-\putLeft) {$u$};
\filldraw[black] (\horDist*\edgeDist,2.5*\edgeDist) circle (2pt) node at (\horDist*\edgeDist,2.5*\edgeDist-\putLeft) {$v$};
\filldraw[black] (2*\edgeDist,2*\edgeDist) circle (0pt) node at (1.5*\edgeDist+\putLeft,1.75*\edgeDist-2.5*\putLeft) {$\mu_{2,1}$};
\filldraw[black] (5*\edgeDist,3*\edgeDist) circle (0pt) node at (3.4*\edgeDist,1.75*\edgeDist-2*\putLeft) {$\mu_{1,2}$};
\filldraw[black] (2*\edgeDist,2*\edgeDist) circle (0pt) node at (1*\edgeDist+\putLeft,3.8*\edgeDist) {$\mu_{1,1}$};
\filldraw[black] (5*\edgeDist,3*\edgeDist) circle (0pt) node at (2.5*\edgeDist+\putLeft,3.9*\edgeDist) {$\mu_{2,2}$};
\filldraw[black] (2*\edgeDist,2*\edgeDist) circle (0pt) node at (1.5*\horDist*\edgeDist,2.5*\edgeDist+.5*\putLeft) {$\omega$};
\filldraw[black] (1.5*\horDist*\edgeDist,2*\edgeDist) circle (0pt) node at (1.5*\horDist*\edgeDist,2.5*\edgeDist-.5*\putLeft) {$\cdots$};
\filldraw[black] (2*\horDist*\edgeDist,2.5*\edgeDist) circle (2pt) node at (2*\horDist*\edgeDist,2.5*\edgeDist+\putLeft) {$r$};
\filldraw[black] (3*\horDist*\edgeDist,0*\edgeDist) circle (2pt) node at (3*\horDist*\edgeDist+\putLeft,0*\edgeDist-\putLeft) {$P_{b}$};
\filldraw[black] (3*\horDist*\edgeDist,1*\edgeDist) circle (2pt) node at (3*\horDist*\edgeDist+\putLeft,1*\edgeDist) {$\vdots$};
\filldraw[black] (3*\horDist*\edgeDist,2*\edgeDist) circle (2pt);
\filldraw[black] (3*\horDist*\edgeDist,3*\edgeDist) circle (2pt);
\filldraw[black] (3*\horDist*\edgeDist,4*\edgeDist) circle (2pt) node at (3*\horDist*\edgeDist+\putLeft,4*\edgeDist) {$\vdots$};
\filldraw[black] (3*\horDist*\edgeDist,5*\edgeDist) circle (2pt) node at (3*\horDist*\edgeDist+\putLeft,5*\edgeDist+\putLeft) {$P_{d}$};
%
%
%
\draw (0*\edgeDist,0*\edgeDist) to (0*\edgeDist,1*\edgeDist); 
\draw (0*\edgeDist,1*\edgeDist) to (0*\edgeDist,2*\edgeDist); 
\draw (0*\edgeDist,3*\edgeDist) to (0*\edgeDist,4*\edgeDist); 
\draw (0*\edgeDist,4*\edgeDist) to (0*\edgeDist,5*\edgeDist); 
\draw (3*\horDist*\edgeDist,0*\edgeDist) to (3*\horDist*\edgeDist,1*\edgeDist); 
\draw (3*\horDist*\edgeDist,1*\edgeDist) to (3*\horDist*\edgeDist,2*\edgeDist); 
\draw (3*\horDist*\edgeDist,3*\edgeDist) to (3*\horDist*\edgeDist,4*\edgeDist); 
\draw (3*\horDist*\edgeDist,4*\edgeDist) to (3*\horDist*\edgeDist,5*\edgeDist); 
\draw (0*\edgeDist,1*\edgeDist) to (2*\horDist*\edgeDist,2.5*\edgeDist); 
\draw (0*\edgeDist,4*\edgeDist) to (1*\horDist*\edgeDist,2.5*\edgeDist); 
\draw (3*\horDist*\edgeDist,4*\edgeDist) to (1*\horDist*\edgeDist,2.5*\edgeDist); 
\draw (1*\horDist*\edgeDist,2.5*\edgeDist) to (2*\horDist*\edgeDist,2.5*\edgeDist); 
\draw (3*\horDist*\edgeDist,1*\edgeDist) to (2*\horDist*\edgeDist,2.5*\edgeDist); 
\draw[dashed] (1*\edgeDist,3*\edgeDist) to (2.5*\edgeDist,3*\edgeDist); 
\draw[dashed] (2*\edgeDist,2*\edgeDist) to (3.5*\edgeDist,2*\edgeDist); 
%
\draw[decorate]  (0*\edgeDist,4*\edgeDist) -- (1*\horDist*\edgeDist,2.5*\edgeDist); 
\draw[decorate]  (1*\horDist*\edgeDist,2.5*\edgeDist) -- (3*\horDist*\edgeDist,4*\edgeDist); 
\draw[decorate] (2*\horDist*\edgeDist,2.5*\edgeDist)  -- (0*\edgeDist,1*\edgeDist); 
\draw[decorate] (3*\horDist*\edgeDist,1*\edgeDist)  -- (2*\horDist*\edgeDist,2.5*\edgeDist); 
\end{tikzpicture}
\subcaption{Case 1}
\end{subfigure}
\hfill
\begin{subfigure}{.45\linewidth}
\begin{tikzpicture}[decoration={brace,amplitude=22pt}]
\filldraw[black] (0*\edgeDist,0*\edgeDist) circle (2pt) node at (-\putLeft,0*\edgeDist-\putLeft) {$P_{c}$};
\filldraw[black] (0*\edgeDist,1*\edgeDist) circle (2pt) node at (-\putLeft,1*\edgeDist) {$\vdots$};
\filldraw[black] (0*\edgeDist,2*\edgeDist) circle (2pt);
\filldraw[black] (0*\edgeDist,3*\edgeDist) circle (2pt);
\filldraw[black] (0*\edgeDist,4*\edgeDist) circle (2pt) node at (-\putLeft,4*\edgeDist) {$\vdots$};
\filldraw[black] (0*\edgeDist,5*\edgeDist) circle (2pt) node at (0*\edgeDist-\putLeft,5*\edgeDist+\putLeft) {$P_{a}$};
\filldraw[black] (3.5*\edgeDist,2*\edgeDist) circle (2pt) node at (3.5*\edgeDist+\putLeft,2*\edgeDist) {$s$};
\filldraw[black] (1*\edgeDist,2*\edgeDist) circle (2pt) node at (1*\edgeDist-\putLeft,2*\edgeDist) {$w$};
\filldraw[black] (3.5*\edgeDist,3*\edgeDist) circle (2pt) node at (3.5*\edgeDist+\putLeft,3*\edgeDist) {$t$};
\filldraw[black] (1*\edgeDist,3*\edgeDist) circle (2pt) node at (1*\edgeDist-\putLeft,3*\edgeDist) {$u$};
\filldraw[black] (\horDist*\edgeDist,2.5*\edgeDist) circle (2pt) node at (\horDist*\edgeDist-\putLeft,2.5*\edgeDist) {$v$};
\filldraw[black] (2*\edgeDist,2*\edgeDist) circle (0pt) node at (1.5*\edgeDist-\putLeft,1.75*\edgeDist-2*\putLeft) {$\mu_{2,1}$};
\filldraw[black] (5*\edgeDist,3*\edgeDist) circle (0pt) node at (3.4*\edgeDist,1.75*\edgeDist-2*\putLeft) {$\mu_{1,2}$};
\filldraw[black] (2*\edgeDist,2*\edgeDist) circle (0pt) node at (1*\edgeDist+\putLeft,3.8*\edgeDist) {$\mu_{1,1}$};
\filldraw[black] (5*\edgeDist,3*\edgeDist) circle (0pt) node at (2*\horDist*\edgeDist+\putLeft,3.8*\edgeDist) {$\mu_{2,2}$};
\filldraw[black] (2*\edgeDist,2*\edgeDist) circle (0pt) node at (1.5*\horDist*\edgeDist,2.5*\edgeDist+.5*\putLeft) {$\omega$};
\filldraw[black] (1.5*\horDist*\edgeDist,2*\edgeDist) circle (0pt) node at (1.5*\horDist*\edgeDist,2.5*\edgeDist-.5*\putLeft) {$\cdots$};
\filldraw[black] (2*\horDist*\edgeDist,2.5*\edgeDist) circle (2pt) node at (2*\horDist*\edgeDist+\putLeft,2.5*\edgeDist) {$r$};
\filldraw[black] (3*\horDist*\edgeDist,0*\edgeDist) circle (2pt) node at (3*\horDist*\edgeDist+\putLeft,0*\edgeDist-\putLeft) {$P_{b}$};
\filldraw[black] (3*\horDist*\edgeDist,1*\edgeDist) circle (2pt) node at (3*\horDist*\edgeDist+\putLeft,1*\edgeDist) {$\vdots$};
\filldraw[black] (3*\horDist*\edgeDist,2*\edgeDist) circle (2pt);
\filldraw[black] (3*\horDist*\edgeDist,3*\edgeDist) circle (2pt);
\filldraw[black] (3*\horDist*\edgeDist,4*\edgeDist) circle (2pt) node at (3*\horDist*\edgeDist+\putLeft,4*\edgeDist) {$\vdots$};
\filldraw[black] (3*\horDist*\edgeDist,5*\edgeDist) circle (2pt) node at (3*\horDist*\edgeDist+\putLeft,5*\edgeDist+\putLeft) {$P_{d}$};
%
%
%
\draw (0*\edgeDist,0*\edgeDist) to (0*\edgeDist,1*\edgeDist); 
\draw (0*\edgeDist,1*\edgeDist) to (0*\edgeDist,2*\edgeDist); 
\draw (0*\edgeDist,3*\edgeDist) to (0*\edgeDist,4*\edgeDist); 
\draw (0*\edgeDist,4*\edgeDist) to (0*\edgeDist,5*\edgeDist); 
\draw (3*\horDist*\edgeDist,0*\edgeDist) to (3*\horDist*\edgeDist,1*\edgeDist); 
\draw (3*\horDist*\edgeDist,1*\edgeDist) to (3*\horDist*\edgeDist,2*\edgeDist); 
\draw (3*\horDist*\edgeDist,3*\edgeDist) to (3*\horDist*\edgeDist,4*\edgeDist); 
\draw (3*\horDist*\edgeDist,4*\edgeDist) to (3*\horDist*\edgeDist,5*\edgeDist); 
\draw (0*\edgeDist,1*\edgeDist) to (1*\horDist*\edgeDist,2.5*\edgeDist); 
\draw (0*\edgeDist,4*\edgeDist) to (1*\horDist*\edgeDist,2.5*\edgeDist); 
\draw (3*\horDist*\edgeDist,4*\edgeDist) to (2*\horDist*\edgeDist,2.5*\edgeDist); 
\draw (1*\horDist*\edgeDist,2.5*\edgeDist) to (2*\horDist*\edgeDist,2.5*\edgeDist); 
\draw (3*\horDist*\edgeDist,1*\edgeDist) to (2*\horDist*\edgeDist,2.5*\edgeDist); 
%
\draw[decorate]  (0*\edgeDist,4*\edgeDist) -- (1*\horDist*\edgeDist,2.5*\edgeDist); 
\draw[decorate]  (2*\horDist*\edgeDist,2.5*\edgeDist) -- (3*\horDist*\edgeDist,4*\edgeDist); 
\draw[decorate] (1*\horDist*\edgeDist,2.5*\edgeDist)  -- (0*\edgeDist,1*\edgeDist); 
\draw[decorate] (3*\horDist*\edgeDist,1*\edgeDist)  -- (2*\horDist*\edgeDist,2.5*\edgeDist); 
\end{tikzpicture}
\subcaption{Case 2}
\end{subfigure}
\caption{Illustration for Proposition \ref{prop_intersection_two_strands_and_four_paths_is_strand}}
\label{fig_intersection_two_strands_and_four_paths_is_strand}
\end{figure}

We are now ready to prove Theorem \ref{thm_doip_iff_no_forbidden_path}.

\bproof[Proof of Theorem \ref{thm_doip_iff_no_forbidden_path}] 
$(1) \implies (2)$: If $\ul{P}$ is DOIP, then in particular $K_{P_{\ind}}$ has no directed two cycle. Corollary \ref{cor_KP_ind_2_cycle_iff_Pind_weakly_forbidden_path} implies that $P_{\ind}$ does not contain an internal strand, an internal fork, a double fork, or a complete ladder.

$(2) \implies (1)$: We assume that $P_{\ind}$ does not contain an internal strand, an internal fork, a double fork, or a complete ladder, and we will show that $\ul{P}$ is DOIP. By Theorem \ref{thm_doip_equiv_K_P_acyclic}, this is equivalent to showing that $K_{P_{\ind}}$ is directed acyclic. Suppose by contradiction that $K_{P_{\ind}}$ has a minimal cycle of length $m$, i.e. that $K_{P_{\ind}}$ has no cycle of size smaller than $m$. Proposition \ref{prop_no_1_cycle} and Corollary \ref{cor_KP_ind_2_cycle_iff_Pind_weakly_forbidden_path} imply that $m \geq 3$. Without loss of generality, we may suppose that $(i,i+1)$ are arcs of $K_{P_{\ind}}$ for $1 \leq i \leq m-1$ and that $(m,1)$ is an arc of $K_{P_{\ind}}$. Proposition \ref{prop_arc_iff_strand} implies that there are strands $\gamma_{i}$ from $P_{i}$ to $P_{i+1}$ for $1 \leq i \leq m-1$ and a strand $\gamma_{m}$ from $P_{m}$ to $P_{1}$. Lemma \ref{lem_strand_intersecting_path_creates_new_arc} and Propositions \ref{prop_intersection_two_strands_is_strand} and \ref{prop_intersection_two_strands_and_four_paths_is_strand} imply for $i,j \in [m]$ that:
\begin{enumerate}
    \item $V(\gamma_{i}) \cap V(\gamma_{j}) = \varnothing$ whenever $\card{i-j} \geq 2$,
    \item $V(\gamma_{i}) \cap V(\gamma_{i+1}) = \varnothing$, and 
    \item $V(\gamma_{i}) \cap V(P_{j}) = \varnothing$ whenever $j\neq i$ or $j \neq i+1$.
\end{enumerate}
If not, these propositions and lemmas would imply the existence of an arc $(i,j)$ for some $i,j \in [m]$ distinct, which would contradict the assumption that $K_{P_{\ind}}$ has a minimal cycle of length $m$. We denote by $R_{i}$ the unique subpath of $P_{i}$ having terminal vertices $\phi_{\gamma_{i-1}}(2)$ and $\phi_{\gamma_{i}}(1)$. It follows that 
\begin{align*}
    Q := R_{1} \ast \gamma_{1} \ast R_{2} \ast \gamma_{2} \ast \cdots \ast \gamma_{m-1} \ast R_{m} \ast \gamma_{m}
\end{align*}
is a cycle. Since $G$ is a block graph, the induced graph on $V(Q)$ is a complete graph. Thus, at most one of the $R_{i}$ contains an internal vertex of $P_{i}$; otherwise, the edge connecting two distinct internal vertices would be an internal strand. Since $m \geq 3$, we may assume without loss of generality that $P_{2}$ and $P_{3}$ are each a singleton edge. It follows from the definition of strands that $R_{2} = P_{2}$ and that $R_{3} = P_{3}$. Thus, the induced subgraph on $V(P_{2}) \cup V(P_{3})$ is a complete ladder, a contradiction.
\eproof

\subsection{Another Characterization of DOIP Paths}

We take the time to record an equivalent characterization of Thereom \ref{thm_doip_iff_no_forbidden_path} in terms of forbidden subgraphs. This reformulation will be particularly useful in the subsequent section as additional constraints are placed upon the forbidden subgraphs.

\bdefn
Let $G$ be a block graph, and let $\ul{P} := P_{1},\ldots,P_{\ell}$ be vertex-disjoint induced paths in $G$. Let $H$ be an internal strand, an internal fork, or a double fork. We say that $H$ is \textbf{edge-disjoint} from $\ul{P}$ if $E(H) \cap E(\ul{P}) = \varnothing$.
\edefn 

\bex
\label{ex_non_edge_disjoint_fork}
Let $G$ be the graph depicted in Figure \ref{fig_forbidden_paths}, and let $\ul{P}$ be as defined in Example \ref{ex_forbidden_subgraphs}. Then, the graph on the vertices $\{2,3,5,10,9,8,7\}$ is a double fork of $P_{\ind}$ which is not edge-disjoint from $\ul{P}$ (because $E(H) \cap E(\ul{P}) = \{ \{9,10\} \}$).
\eex 

\bprop 
\label{prop_doip_iff_no_forbidden_path_another_char}
Let $G$ be a block graph, and let $\ul{P} := P_{1},\ldots,P_{\ell}$ be vertex-disjoint induced paths of $G$. The following are equivalent:
\begin{enumerate}
    \item $P_{\ind}$ does not contain any of the following as subgraphs:
    \begin{enumerate}
        \item an internal strand, 
        \item an internal fork, 
        \item a double fork, 
        \item a complete ladder.
    \end{enumerate}
   \item $P_{\ind}$ does not contain any of the following as subgraphs:
   \begin{enumerate}
       \item an internal strand which is edge-disjoint from $\ul{P}$, 
       \item an internal fork which is edge-disjoint from $\ul{P}$, 
       \item a double fork which is edge-disjoint from $\ul{P}$, 
       \item a complete ladder.
   \end{enumerate}
\end{enumerate}
\eprop 

\bproof
It is clear that (1) implies (2). That (2) implies (1) follows from Lemmas \ref{lem_internal_strand_implies_edge_disjoint_internal_strand} and \ref{lem_internal_fork_implies_edge_disjoint_internal_fork}.
\eproof

\blem
\label{lem_internal_strand_implies_edge_disjoint_internal_strand}
If $P_{\ind}$ contains an internal strand as a subgraph, then $P_{\ind}$ contains an internal strand that is edge-disjoint from $\ul{P}$.
\elem

\bproof
Let $H$ be an internal strand of $P_{\ind}$. We prove by induction on $\card{E(H)}$ that there exists an internal strand of $P_{\ind}$ that is edge-disjoint from $\ul{P}$. If $\card{E(H)} = 1$, then by Definition \ref{defn_strand}, $H$ is edge-disjoint from $\ul{P}$. Suppose that $H$ is an internal strand with $\card{E(H)} \geq 2$. If $H$ is edge-disjoint from $\ul{P}$, there is nothing to prove. Suppose that $E(P_{i}) \cap E(H) \neq \varnothing$ for some $1 \leq i \leq \ell$. Lemma \ref{lem_decomp_intersecting_induced_paths} implies that $H = H_{1} \ast \gamma \ast H_{2}$ where $\gamma$ is a subpath of $P_{i}$ with $\card{E(\gamma)} \geq 1$. Because $H$ does not contain $P_{i}$, at least one of the terminal vertices of $H_{1}$ or $H_{2}$ is an internal vertex of $P_{i}$. Without loss of generality, suppose that it is $H_{1}$. Then, $H_{1}$ is an internal strand of $P_{\ind}$ with $\card{E(H_{1})} < \card{E(H)}$. By the induction hypothesis, applied to $H_{1}$, there exists an internal strand of $P_{\ind}$ which is edge-disjoint from $\ul{P}$, which completes the proof.
\eproof 

\blem
\label{lem_internal_fork_implies_edge_disjoint_internal_fork}
If $P_{\ind}$ contains an internal fork or a double fork as a subgraph, then $P_{\ind}$ contains an internal strand, an internal fork, or a double fork which is edge-disjoint from $\ul{P}$.
\elem

\bproof
Let $H$ be an internal fork or a double fork of $P_{\ind}$. Let $Q$ be the subpath of $H$, as defined in Definitions \ref{defn_fork} and \ref{defn_double_fork}. The proof proceeds by induction on $\card{E(Q)}$. When $\card{E(Q)} = 0$, i.e., $Q$ is the path consisting of a singleton vertex, it is clear from the definitions that $H$ is edge-disjoint from $\ul{P}$. When $\card{E(Q)} \geq 1$, the proof proceeds via induction, as in the proof of Lemma \ref{lem_internal_strand_implies_edge_disjoint_internal_strand}. We observe that in replicating the proof of Lemma \ref{lem_internal_strand_implies_edge_disjoint_internal_strand} that: if $H$ is an internal fork, then $H_{1}$ or $H_{2}$ is a strand or a fork with at least one of them being internal; and if $H$ is a double fork, then $H_{1}$ and $H_{2}$ are forks with at least one of them being internal.
\eproof

\bex
\label{ex_non_edge_disjoint_fork_splits_two_forks}
We illustrate Lemma \ref{lem_internal_fork_implies_edge_disjoint_internal_fork} in the context of Example \ref{ex_non_edge_disjoint_fork}. Observe that $H_{1}$ and $H_{2}$ are the graphs on the vertices $\{2,3,5,10\}$ and $\{9,8,7\}$, respectively. Then, $H_{1}$ is an internal fork of $P_{\ind}$.
\eex 

\section{Combinatorial Characterization of $\reg(S/J_{G})$ for Block Graphs}
\label{sec_combinatorial_char_reg_for_block_graphs}

In this section, we prove the following theorem. 

\bthm 
\label{thm_reg=nu_block_graph}
Let $G$ be a block graph. Then,
\begin{align*}
    \nu(G) = \reg(S/J_{G}).
\end{align*}
\ethm 

It suffices to show that $\reg(S/J_{G}) \leq \nu(G)$, and to do so, we utilize the theory developed in \cite{malayeri2021proof}, which we now recall.

\bdefn
\label{defn_relevant_compatibility}
For a graph $G$, define $\wh{G} := G \setminus I_{S}(G)$ where $I_{S}(G)$ denotes the set of isolated vertices of $G$. For $v \in V(G)$, recall that $N_{G}(v)$ denotes the vertices of $G$ adjacent to $v$. For $v \in V(G)$, we define the graph $G_{v}$ as follows:
\begin{align*}
V(G_{v}) &:= V(G) \\
E(G_{v}) &:= E(G) \cup \{ \{u,w\} \mid u,w \in N_{G}(v)\}.
\end{align*}
The graph $G_{v}$ is referred to as the \textbf{completion of $G$ at $v$}.
\edefn

\bdefn
We will say that a subset $\G$ of all finite graphs is \textbf{compatible} if $\G$ satisfies the following conditions:
\begin{enumerate}
\item $\bigsqcup_{i=1}^{t} K_{n_{i}} \in \G$ for all $n_{i} \in \Z$ with $n_{i} \geq 2$,
\item $\wh{G} \in \G$ for all $G \in \G$,
\item $G \setminus \{v\} \in \G$ for all $G \in \G$ and $v \in V(G)$,
\item $G_{v} \in \G$ for all $G \in \G$ and $v \in V(G)$.
\end{enumerate}
\edefn

\bdefn[{\cite[Definition 2.1]{malayeri2021proof}}]
\label{defn_comp_map}
Let $\G$ be a subset of all finite graphs. Suppose that $\G$ is compatible. A map $\varphi : \G \ra \N_{0}$ is called \textbf{compatible} if it satisfies the following conditions:
\begin{enumerate}
\item $\varphi(\wh{G}) \leq \varphi(G)$ for all $G \in \G$, \label{defn_comp_map_1}
\item if $G = \bigsqcup_{i=1}^{t} K_{n_{i}}$, where $n_{i} \geq 2$ for every $1 \leq i \leq t$, then $\varphi(G) \geq t$, \label{defn_comp_map_2}
\item if $G \neq \bigsqcup_{i=1}^{t} K_{n_{i}}$, then there exists $v \in V(G)$ such that 
\begin{enumerate}
\item $\varphi(G \setminus v) \leq \varphi(G)$, and \label{defn_comp_map_3a}
\item $\varphi(G_{v}) < \varphi(G)$. \label{defn_comp_map_3b}
\end{enumerate}
\end{enumerate}
\edefn

\bthm[{\cite[Theorem 2.3]{malayeri2021proof}}]
\label{thm_reg_leq_phi}
Let $\G$ be a subset of all finite graphs. Suppose that $\G$ is compatible and that $\varphi : \G \ra \N_{\geq 0}$ is compatible. Then, for all $G \in \G$,
\begin{align*}
\reg(S/J_{G}) \leq \varphi(G).
\end{align*} 
\ethm

\brem
In \cite{malayeri2021proof}, they proved \cref{thm_reg_leq_phi} where $\G$ is the set of all finite graphs. Their proof technique has two main steps. First, they show for a graph $G$ that 
\begin{align*}
\reg(S/J_{G\setminus \{v\}}) &\leq \varphi(G \setminus \{v\}) \\
\reg(S/J_{G_{v}}) &< \varphi(G_{v})
\end{align*}
Second, they utilize induction on the number of internal vertices of a graph together with the short exact sequence
\begin{align*}
    0 \ra S/J_{G} \ra S/J_{G_{v}} \oplus S_{v}/J_{G \setminus v} \ra S_{v}/J_{G_{v} \setminus v} \ra 0
\end{align*}
to deduce that $\reg(S/J_{G}) \leq \varphi(G)$. The fact that $\mathcal{G}$ is compatible, i.e., closed under vertex completion and deletion, allows us to apply the induction step in our setting.
\erem

\bex
\label{ex_chor_graphs_comp}
The class of chordal graphs and the class of block graphs are compatible.
\eex

\blem
\label{lem_nu_satisfies_1_2_3a}
Let $\G$ be a compatible subset of finite graphs, and let $\nu : \G \ra \N_{0}$ be defined as in Definition \ref{defn_nu}. Then, $\nu$ satisfies conditions \ref{defn_comp_map_1}, \ref{defn_comp_map_2}, and \ref{defn_comp_map_3a} of Definition \ref{defn_comp_map}.
\elem

\bproof
The claims follow from the straightforward observations that:
\begin{enumerate}
    \item If $H$ is an induced subgraph of $G$, then $\nu(H) \leq \nu(G)$,
    \item If $G = G_{1} \sqcup G_{2}$, then $\nu(G) = \nu(G_{1}) + \nu(G_{2})$,
    \item $\nu(K_{n}) = 1$ for all $n \geq 2$.
\end{enumerate}
\eproof

Theorem \ref{thm_reg_leq_phi} and Lemma \ref{lem_nu_satisfies_1_2_3a} show that to prove Theorem \ref{thm_reg=nu_block_graph} it suffices to show that there exists a vertex $c$ of the block graph $G$ such that $\nu(G_{c}) < \nu(G)$. We now introduce some notation and prove a few preparatory lemmas.

\bdefn
\label{defn_two_block_prop}
A vertex $v$ of $G$ is called a \textbf{cut vertex} if the number of connected components of $G \setminus \{v\}$ is strictly larger than the number of connected components of $G$. For a vertex $v$ of $G$, we define the \textbf{clique degree of $v$}, denoted by $\cdeg(v)$, as the number of maximal distinct cliques of $G$ containing $v$. The number of maximal cliques in a block graph $G$ is denoted by $c(G)$. We say that a block graph $G$ has the \textbf{two-block property} if, for every vertex $v$ of $G$, $\cdeg(v) \leq 2$. We say that a block graph $G$ is a \textbf{path of cliques} if every block of $G$ has at most two cut vertices. In a path of cliques, blocks with exactly one cut vertex are called \textbf{terminal blocks}.
\edefn

\bex 
Consider the graph constructed from the complete graph on three vertices by attaching a whisker to each vertex of the complete graph. (This graph is sometimes referred to as the \textbf{net}.) This graph has the two-block property but is not a path of cliques.
\eex 

\blem
\label{lem_all_edges_admissible_two_block_graphs}
Let $G$ be a block graph that has the two-block property. Let $\A$ denote a collection of edges of $G$ such that no two edges belong to a common clique of $G$. Let $\ul{P}$ denote the disjoint union of paths obtained from $\mathcal{A}$ after concatenating those edges of $\mathcal{A}$ sharing a terminal vertex. Let $P_{\ind}$ denote the induced subgraph of $G$ on $V(\ul{P})$. Then, $P_{\ind}$ is DOIP.
\elem

\bproof
Suppose by contradiction that there exists a block graph $G$ with the two-block property and vertex-disjoint edges $\mc{A}$ such that $P_{\ind}$ is not DOIP. We may suppose that among all such block graphs that $G$ has been chosen to minimize $c(G)$. It is clear that $c(G) \geq 3$, as any block graph on two or fewer cliques, together with any choice of edges $\A$, realizes $P_{\ind}$ that is DOIP. By Theorem \ref{thm_doip_iff_no_forbidden_path} and Proposition \ref{prop_doip_iff_no_forbidden_path_another_char}, $P_{\ind}$ contains $H$, an internal strand, an internal fork, or a double fork, which is edge-disjoint from $\ul{P}$ ($H$ is not a complete ladder because no two edges of $\A$ belong to the same clique). By minimality of $c(G)$, we may assume that 
\begin{enumerate}
    \item $V(H) \cap V(B_{i}) \neq \varnothing$ for all $1 \leq i \leq c(G)$, and
    \item for every block $B$ of $G$, there exists an edge $e$ of $\A$ contained in $B$.
\end{enumerate}
(If $G$ did not satisfy these two conditions, then we could produce a smaller counterexample by deleting irrelevant blocks of $G$.) The first condition implies that $G$ is a path of cliques, since $H$ is necessarily contained in a path of cliques. We may label the cliques of $G$ consecutively, starting from a terminal clique, by $B_{1},B_{2},B_{3},\ldots,B_{c(G)}$. We denote by $v_{i}$ the vertex common to $B_{i}$ and $B_{i+1}$ for $1 \leq i \leq c(G)-1$. We denote by $e_{i}$ the edge of $\A$ belonging to $B_{i}$. Consequently, $H$ contains the subpath $[v_{1},v_{2},v_{3},\ldots,v_{c(G)-1}]$. In particular, $v_{i} \in V(P_{\ind})$. We will show by considering two cases below that $v_{2}$ is not a vertex of $e_{2}$. In which case, we can construct $G^{'}$ by deleting the block $B_{1}$ from $G$, $\A^{'}$ by deleting $e_{1}$ from $\A$, and $\ul{P}^{'}$ and $P_{\ind}^{'}$ coming from $\A^{'}$ and $G^{'}$. If $e_{2} = \{a,b\}$, then we construct $H^{'}$ by:
\begin{align*}
    V(H^{'}) &:= \left( V(H) \setminus \{v_{1}\} \right) \cup \{a,b\} \\
    E(H^{'}) &:= \left( E(H \setminus \{v_{1}\}) \right) \cup \{ \{v_{2},a\}, \{v_{2},b\} \}
\end{align*}
is an internal fork or a double fork of $P_{\ind}^{'}$. This would contradict minimality of $c(G)$.

\ul{Claim}: $e_{2}$ does not contain $v_{2}$.

\ul{Proof of Claim}. We consider the following two cases.

\ul{Case 1.} Suppose that $e_{1}$ does not contain $v_{1}$. Then, it must be the case that $e_{2}$ contains $v_{1}$ (because $V(H) \subset V(\ul{P})$). As $e_{2} \neq \{v_{1},v_{2}\}$ (because $H$ is edge-disjoint from $\ul{P}$), $e_{2}$ does not contain $v_{2}$.

\ul{Case 2.} Suppose that $e_{1}$ contains $v_{1}$. Then, $E(H) \cap E(B_{1}) = \varnothing$, since $H$ does contain $e_{1}$ and $G$ is a path of cliques having $B_{1}$ as a terminal vertex. It follows that $v_{1}$ is a terminal vertex of $H$. Hence, it must be the case that $e_{2}$ contains $v_{1}$. Since $e_{2} \neq \{v_{1},v_{2}\}$, it again follows that $e_{2}$ does not contain $v_{2}$.
\eproof

\bcor
\label{cor_nu=reg_2_block_prop}
If $G$ is a block graph with the two-block property, then $\nu(G) = \reg(S/J_{G}) = c(G)$.
\ecor 

\brem
In Corollary \ref{cor_nu=reg_2_block_prop}, the statement that $\reg(S/J_{G}) = c(G)$ follows from {\cite[Proposition 1.3]{herzog2018extremal}}.
\erem

\bproof
For this family of graphs, Lemma \ref{lem_all_edges_admissible_two_block_graphs} proves that that $c(G) \leq \nu(G)$. For any graph $G$, we have that $\nu(G) \leq \reg(S/J_{G}) \leq c(G)$ by Theorem \ref{thm_nu_leq_reg} and {\cite[Corollary 2.7]{malayeri2021proof}}.
\eproof

\blem 
\label{lem_terminal_vertices_of_path_of_cliques}
Let $G$ be a block graph. Suppose that $G = G_{1} \cup G_{2}$ and that $G_{1} \cap G_{2} = \{c\}$ for some vertex $c$ of $G$ where $G_{1}$ has the two-block property and $G_{2}$ is a block graph. Let $v_{1} \in V(G_{1})$. Suppose that $\ul{P}$ is a union of vertex-disjoint induced paths of $G$ that contains the edge $\{v_{1},c\}$. If $H$ is an internal strand, an internal fork, or a double fork of $P_{\ind}$ that is edge-disjoint from $\ul{P}$, then 
\begin{align*}
    V(H) \cap (V(G_{1})\setminus \{c\}) = \varnothing.
\end{align*}
\elem 

\bproof
Suppose by contradiction that $V(H) \cap (V(G_{1})\setminus \{c\}) \neq \varnothing$. Then, $c$ is an internal vertex of $H$; otherwise, $H$ would be a subgraph of $G_{1}$, which would contradict Lemma \ref{lem_all_edges_admissible_two_block_graphs} together with Theorem \ref{thm_doip_iff_no_forbidden_path}. We denote by $a$ the vertex of $V(H) \cap V(G_{1})$ that is adjacent to $c$. Since $\{v_{1},c\} \notin E(H)$, it must be the case that $a \neq v_{1}$. We define the graph $H^{'}$ by 
\begin{align*}
    V(H^{'}) &:= \left( V(H) \cap V(G_{1}) \right) \cup \{v_{1}\} \\
    E(H^{'}) &:= \left( E(H) \cap E(G_{1}) \right) \cup \{ \{a,v_{1}\} \}.
\end{align*}
It follows from Definitions \ref{defn_fork}, \ref{defn_double_fork} that $H^{'}$ is an internal fork or a double fork of $P_{\ind}$. This contradicts Lemma \ref{lem_all_edges_admissible_two_block_graphs}. 
\eproof 

\blem 
\label{lem_reduction_to_c_terminal_vertex}
Let $G$ be a block graph, $\ul{P}$ vertex-disjoint induced paths of $G$, and $H$ an internal strand, an internal fork, or a double fork of $P_{\ind}$ that is edge-disjoint from $\ul{P}$. Suppose that $c$ is an internal vertex of both $H$ and $\ul{P}$. Then, there exists an internal strand or an internal fork of $P_{\ind}$ which is edge-disjoint from $\ul{P}$ and which contains $c$ as a terminal vertex.
\elem 

\bproof
Since $c$ is an internal vertex of $H$, there exist subgraphs $H_{1}$ and $H_{2}$ of $H$ such that 
\begin{align*}
    H &= H_{1} \cup H_{2} \\
    \{c\} &= H_{1} \cap H_{2}.
\end{align*}
Now, $H_{1}$ is an internal strand or an internal fork of $P_{\ind}$ which is edge-disjoint from $\ul{P}$.
\eproof

\bprop
\label{prop_nu_decreases_along_somme_completions}
Let $G$ be a block graph. Then, for some cut vertex $c$ of $G$, we have that $\nu(G_{c}) < \nu(G)$.
\eprop 

\bproof
If $G$ has the two-block property, then the result follows from Corollary \ref{cor_nu=reg_2_block_prop}. Thus, we may assume that $G$ is a block graph which does not have the two-block property. Pick $c \in \Cut(G)$ such that $\cdeg(c) \geq 3$ and $G_{1},G_{2},\ldots,G_{t}$ are subgraphs of $G$, $t \geq 3$, satisfying:
\begin{enumerate}
    \item $G_{i}$ has the two-block property for $1 \leq i \leq t-1$,
    \item $G = \bigcup_{i=1}^{t} G_{i}$, and
    \item $G_{i} \cap G_{j} = \{c\}$ for $1 \leq i < j \leq t$. 
\end{enumerate}
Such a $c$ exists because $G$ does not have the two-block property, and by induction on the number of blocks of $G$. For $1 \leq i \leq t$, denote by $B_{i}$ the block of $G_{i}$ that contains $c$. Let $B$ be the block of $G_{c}$ containing $c$. Let $\ul{P}^{'}$ be DOIP paths of $G_{c}$ such that 
\begin{align*}
    \nu(G_{c}) = \card{E(\ul{P}^{'})}.
\end{align*}
Let $P_{\ind}^{'}$ be the induced subgraph of $G_{c}$ on $V(\ul{P}^{'})$. From $\ul{P}^{'}$, we will construct a DOIP path $\ul{P}$ of $G$ such that 
\begin{align*}
    \card{E(\ul{P}^{'})} < \card{E(\ul{P})}.
\end{align*}
From which, it will follow that $\nu(G_{c}) < \nu(G)$. This construction will proceed across several cases.

\vspace{1em}

\ul{Case 1.} For this case, refer to Figure~\ref{fig_prop_nu_decreases_along_somme_completions}. Suppose that $E(B) \cap E(\ul{P}^{'}) = \varnothing$. It follows that $\ul{P}^{'}$ is a subgraph of $G$. In particular, it follows that $E(B_{i}) \cap E(\ul{P}^{'}) = \varnothing$ for all $1 \leq i \leq t$. Pick $v_{1} \in B_{1} \setminus \{c\}$. We define $\ul{P}$ to be the subgraph of $G$ as follows: 
\begin{align*}
V(\ul{P}) &:= V(\ul{P}^{'}) \cup \{v_{1},c\} \\
E(\ul{P}) &:= E(\ul{P}^{'}) \cup \{ \{v_{1},c\} \}.
\end{align*}
The assumptions that $G_{1}$ has the two-block property and that $E(B_{1}) \cap E(\ul{P}^{'}) = \varnothing$ imply that $\ul{P}$ consists of vertex-disjoint induced paths. Suppose by contradiction that $H$ is an internal strand, an internal fork, or a double fork of $P_{\ind}$ which is edge-disjoint from $E(\ul{P})$. Lemma \ref{lem_terminal_vertices_of_path_of_cliques} implies that $H$ does not contain $v_{1}$. Hence, $H$ contains $c$; otherwise, $H$ would be an internal strand, an internal fork, or a double fork of $P_{\ind}^{'}$, a contradiction. Moreover, $c$ is not a terminal vertex of $H$, since $\deg_{\ul{P}}(c) = 1$ and $V(H) \cap (V(B_{1}) \setminus \{c\}) = \varnothing$. Hence, $c$ is an internal vertex of $H$. We denote the vertices of $H$ adjacent to $c$ by $a_{1}$ and $a_{2}$, and we define the graph $H^{'}$ of $G_{c}$ as follows:
\begin{align*}
    V(H^{'}) &:= V(H) \setminus \{c\} \\
    E(H^{'}) &:= \left( E(H) \setminus \{ \{a_{1},c\}, \{a_{2},c\} \} \right) \cup \{ \{a_{1},a_{2}\} \}.
\end{align*}
Because $a_{1}$ and $a_{2}$ are adjacent to $c$, $\{a_{1},a_{2}\}$ is indeed an edge of $G_{c}$. We observe that the condition $E(B) \cap E(\ul{P}^{'}) = \varnothing$ implies that $\{a_{1},a_{2}\}$ does not contain an edge of $\ul{P}^{'}$. Moreover, this condition implies that if $a_{i}$ is a terminal vertex of $H$, then $a_{i}$ is an internal vertex of $\ul{P}$. It now follows from Definitions \ref{defn_strand}, \ref{defn_fork}, \ref{defn_double_fork} that $H^{'}$ is an internal strand, an internal fork, or a double fork of $P_{\ind}^{'}$, a contradiction.

\begin{figure}
\begin{subfigure}{.45\linewidth}
\centering
\begin{tikzpicture}
\pgfmathsetmacro\R{2.5}
\pgfmathsetmacro\putLeft{.5*.8}
    \draw[smooth,black] plot[domain=0:36*5,samples=200] (\x:{\R*cos(3*\x)});
    \filldraw[black] (0,0) circle (2pt) node at (-.5,0) {$c$};
    \filldraw[black] (\R,0) circle (2pt) node at (\R+.5,0) {$a_{2}$};
    \filldraw[black] (\R*-0.5,\R*0.866) circle (2pt) node at (\R*-0.5 -.3,\R*0.866 + .3) {$v_{1}$};
    \filldraw[black] (\R*-0.5,\R*-.866) circle (2pt) node at (\R*-0.5 -0.3,\R*-.866 -0.3) {$a_{1}$};
    \draw (0,0) to (\R*-0.5,\R*0.866); 
    \draw[dashed] (0,0) to (\R*-0.5,\R*-.866); 
    \draw[dashed] (0,0) to (\R,0); 
\end{tikzpicture}
\subcaption{$G$, $P_{\ind}$ (solid line), and $H$ (dashed lines)}
\end{subfigure}
\hfill 
\begin{subfigure}{.45\linewidth}
\centering
\begin{tikzpicture}
\pgfmathsetmacro\R{2}
    \draw[black] (0,0) circle (\R);
    \filldraw[black] (0,0) circle (2pt) node at (-.4,0) {$c$};
    \filldraw[black] (\R,0) circle (2pt) node at (\R+.5,0) {$a_{2}$};
    \filldraw[black] (\R*-0.5,\R*0.866) circle (2pt) node at (\R*-0.5 -.3,\R*0.866 + .3) {$v_{1}$};
    \filldraw[black] (\R*-0.5,\R*-.866) circle (2pt) node at (\R*-0.5 -0.3,\R*-.866 -0.3) {$a_{1}$};
    \draw[dashed] (\R*-0.5,\R*-.866) to (\R,0); 
\end{tikzpicture}
\subcaption{$G_{c}$ and $H^{'}$ (dashed line)}
\end{subfigure}
\caption{Illustration for Case 1}
\label{fig_prop_nu_decreases_along_somme_completions}
\end{figure}
We next consider the case where $E(B) \cap E(\ul{P}^{'}) \neq \varnothing$. We distinguish between cases based on whether this edge of $E(B) \cap E(\ul{P}^{'})$ contains the vertex $c$.

\ul{Case 2.} Suppose that $\{a,c\}$ is an edge of $E(B) \cap E(\ul{P}^{'})$. If necessary, relabel $B_{1}$ and $G_{1}$ by $B_{2}$ and $G_{2}$, respectively, so that we may assume that $a \notin V(G_{1})$. We pick $v_{1} \in V(B_{1}) \setminus \{c\}$ and construct $\ul{P}$ of $G$ as follows:
\begin{align*}
    V(\ul{P}) &:= V(\ul{P}^{'}) \cup \{v_{1}\} \\
    E(\ul{P}) &:= E(\ul{P^{'}}) \cup \{ \{v_{1},c\} \}.
\end{align*}
We observe that $\ul{P}$ consists of vertex-disjoint induced paths of $G$, since $E(\ul{P}) \cap E(B_{1}) = \varnothing$ and $G_{1}$ has the two-block property. 
Suppose by contradiction that $H$ is an internal strand, an internal fork, or a double fork of $P_{\ind}$ which is edge-disjoint from $E(\ul{P})$. Lemma \ref{lem_terminal_vertices_of_path_of_cliques} implies that $H$ does not contain $v_{1}$. Hence, $H$ contains $c$; otherwise, $H$ would be an internal strand, an internal fork, or a double fork of $P_{\ind}^{'}$, a contradiction. By Lemma \ref{lem_reduction_to_c_terminal_vertex}, we may assume that $c$ is a terminal vertex of $H$. Let $b$ denote the vertex of $H$ which is adjacent to $c$. We construct $H^{'}$ of $G_{c}$ as follows:
\begin{align*}
    V(H^{'}) &:= V(H) \cup \{b\} \\
    E(H^{'}) &:= E(H) \cup \{ \{b,a\} \}.
\end{align*}
Since $H$ is edge-disjoint from $\ul{P}$, $b \neq a$. Hence, $H^{'}$ is an internal fork or a double fork of $P_{\ind}^{'}$, a contradiction.

\ul{Case 3.} Suppose that $\{a,b\}$ is an edge of $E(B) \cap E(\ul{P}^{'})$ with $a \neq c$ and $b \neq c$. Let $v_{1}$ and $v_{2}$ be vertices of $B_{1} \setminus \{c\}$ and $B_{2} \setminus \{c\}$, respectively. We construct $\ul{P}$, vertex-disjoint induced paths of $G$ from $\ul{P}^{'}$, by deleting the edge $\{a,b\}$ from $\ul{P}^{'}$, removing any isolated vertices created after deleting this edge, and then adding the edges $\{v_{1},c\}$ and $\{v_{2},c\}$. Suppose by contradiction that $H$ is an internal strand, an internal fork, or a double fork of $P_{\ind}$ which is edge-disjoint from $\ul{P}$. Lemma \ref{lem_terminal_vertices_of_path_of_cliques} allows us to assume that $H$ contains $c$, and Lemma \ref{lem_reduction_to_c_terminal_vertex} allows us to assume that $c$ is a terminal vertex of $H$. We denote by $d$ the vertex of $H$ which is adjacent to $c$.

\ul{Subcase (a)}. We suppose that $d \neq a$ and that $d \neq b$. Then, we construct $H^{'}$, an internal fork or a double fork of $P_{\ind}^{'}$, as follows:
\begin{align*}
    V(H^{'}) &:= (V(H) \setminus \{c\}) \cup \{a,b\} \\
    E(H^{'}) &:= (E(H) \setminus \{ \{c,d\} \}) \cup \{ \{a,d\}, \{b,d\} \}.
\end{align*}

\ul{Subcase (b)}. We suppose without loss of generality that $a = d$. This implies that $a \in V(\ul{P})$. The construction of $\ul{P}$ from $\ul{P}^{'}$ involved deleting the edge $\{a,b\}$ and any isolated vertices. Hence, it must be the case that $\deg_{\ul{P}^{'}}(a) = 2$, i.e., that $a$ is an internal vertex of $\ul{P}^{'}$. We construct the graph $H^{'}$ as follows:
\begin{align*}
    V(H^{'}) &:= V(H) \setminus \{c\} \\
    E(H^{'}) &:= E(H) \setminus \{ \{a,c\} \}.
\end{align*}
We observe that $H$ being edge-disjoint from $\ul{P}$ implies that if $a \in V(P_{i})$, then $V(H^{'}) \cap V(P_{i}) = \{a\}$. Hence, $H^{'}$ is an internal strand or an internal fork of $P_{\ind}^{'}$, a contradiction.
\eproof


\section{Acknowledgements}

The author would like to sincerely thank the referees for their thorough feedback, which greatly improved the quality of this paper. The author would like to thank his advisor, Uli Walther, for helpful feedback on an earlier draft of this paper.

\bibliographystyle{amsalpha}
\bibliography{bibliography_paper_regularity_binomial_edge_ideals}

\end{document}